\pdfoutput=1
\documentclass[a4paper,11pt]{mynotes}

\usepackage{geometry}
\usepackage{graphicx}
\usepackage{amsfonts}
\usepackage{fontawesome}
\usetikzlibrary{patterns,hobby,calc,intersections}
\usepgfplotslibrary{fillbetween}

\newbool{proofInMain}
\newbool{proofInAppendix}
\newbool{duplicateStatementInAppendix}
\newbool{showLinks}
\newbool{showBoxes}

\boolfalse{proofInMain}
\booltrue{proofInAppendix}
\boolfalse{duplicateStatementInAppendix}
\boolfalse{showLinks}
\boolfalse{showBoxes}

\let\gradient\relax

\let\tangent\relax

\let\norm\relax

\newcommand{\naturals}                       {\mathbb{N}}

\newcommand{\reals}                          {\mathbb{R}}
\newcommand{\realsBar}                       {\bar{\reals}}
\newcommand{\nonnegativeReals}               {\reals_{\geq 0}}

\renewcommand{\d}                            {\mathrm{d}}

\newcommand{\st}                             {\,|\,}

\DeclareMathOperator{\epi}                  {epi}
\newcommand{\epigraph}[1]                   {\epi#1}

\newcommand{\innerProduct}[2]                {\left\langle #1,#2\right\rangle}
\newcommand{\norm}[1]                        {\left\Vert #1\right\Vert}

\newcommand{\transpose}[1]                   {{#1}^\top}

\newcommand{\Pp}[2]                          {\mathcal{P}_{#1}(#2)}
\newcommand{\PpLocal}[3]                     {\mathcal{P}_{#1}(#2)_{#3}}
\newcommand{\Ppabs}[2]                       {\mathcal{P}_{#1,\mathrm{abs}}(#2)}

\DeclareMathOperator{\supp}                  {supp}
\newcommand{\pushforward}[2]                {{#1}_{\#}{#2}}
\newcommand{\gaussian}[2]                {\mathcal{N}\left(#1, #2\right)}

\newcommand{\almostEverywhere}[1]            {#1\text{-a.e.}}

\newcommand{\expectedValue}[2]               {\mathbb{E}^{#1}\ifstrempty{#2}{}{\left[#2\right]}}
\newcommand{\variance}[2]                    {\mathrm{Var}^{#1}\left[#2\right]}
\DeclareMathOperator{\kullbackLeibler}{KL}
\newcommand{\kl}[3]             {\kullbackLeibler\ifstrempty{#3}{}{[{#3}]}(#1||#2)}

\newcommand{\wassersteinDistance}[3]         {W_{#3}\ifstrempty{#1}{}{\left(#1, #2\right)}}
\newcommand{\setPlans}[2]                    {\Gamma(#1\ifstrempty{#2}{}{,#2})}
\newcommand{\setPlansOptimalCost}[3]         {\Gamma_{#3}(#1,#2)}
\newcommand{\setPlansOptimal}[2]             {\setPlansOptimalCost{#1}{#2}{o}}
\newcommand{\setPlansCommonMarginal}[3]      {\Gamma^{#3}(#1,#2)}
\newcommand{\wassersteinBall}[3]             {\mathbb{B}_{\wassersteinDistance{}{}{#3}}(#1;#2)}
\newcommand{\closedWassersteinBall}[3]       {\bar{\mathbb{B}}_{\wassersteinDistance{}{}{#3}}(#1;#2)}

\newcommand{\projection}[2]                  {\pi_{#1}\ifstrempty{#2}{}{\left[#2\right]}}
\newcommand{\comp}[2]                        {#2 \circ #1}
\newcommand{\identity}[1]                    {\mathrm{Id}_{#1}}

\newcommand{\convergenceNarrow}[0]           {\rightharpoonup}
\newcommand{\convergenceWasserstein}[0]           {\to}

\newcommand{\fromAbove}[0]                   {\searrow}

\newcommand{\closure}[2]                {\overline{#1}^{#2}}

\newcommand{\boundary}[2]                    {\partial_{#2} #1}

\newcommand{\Ccinf}[1]                       {C_c^\infty(#1)}
\newcommand{\Cb}[1]                          {C_b\ifstrempty{#1}{}{(#1)}}
\newcommand{\Lp}[3]                          {L^{#1}(#2;#3)}
\newcommand{\C}[3]                           {C^{#1}\ifstrempty{#2}{}{\ifstrempty{#3}{(#2)}{(#2,#3)}}}
\newcommand{\gradient}[2]                   {\nabla_{#1}{#2}}
\newcommand{\gradientStrong}[2]             {\nabla^{s}_{#1}{#2}}

\DeclareMathOperator{\tangent}               {T}

\newcommand{\tangentCone}[2]                 {{\tangent_{#1}(#2)}}
\newcommand{\tangentSpaceInterior}[2]        {\tangent^{\mathrm{int}}_{#1}(#2)}
\DeclareMathOperator{\normal}                {N}
\newcommand{\normalConeRegular}[2]           {\hat{\normal}_{#1}\ifstrempty{#2}{}{(#2)}}

\newcommand{\normalCone}[2]                  {\normal_{#1}\ifstrempty{#2}{}{(#2)}}
\newcommand{\normalConeStrong}[2]            {\normal_{#1}^s\ifstrempty{#2}{}{(#2)}}
\newcommand{\subgradientRegular}[1]          {\hat{\partial}#1}
\newcommand{\subgradientRegularStrong}[1]    {\hat{\partial}^s#1}
\newcommand{\subgradient}[1]                 {\partial#1}

\newcommand{\subgradientHzn}[1]              {\partial^\infty#1}

\newcommand{\WikipediaLink}[2]{\ifbool{showLinks}{\href{#1}{\color{ForestGreen}{#2}}}{#2}}%
\newcommand{\StackExchangeLink}[2]{\ifbool{showLinks}{\href{#1}{\color{purple}{#2}}}{#2}}%

\newcommand{\includeElement}[4]{
\ifbool{#3}{
\begin{#1}[\input{elements/#2_title}]
\uniqueLabel{#1:#2}
\input{elements/#2_body}
\end{#1}
}{}
\ifbool{#4}{
\begin{proof}
{\ifbool{#3}{Proof}{Proof of \cref{#1:#2}
}.}
\input{elements/#2_proof}
\end{proof}
}{}
}
\newcommand{\uniqueLabel}[1]{\label{#1}}

\newcommand{\emphBox}[2]{%
\ifbool{showBoxes}{%
\begin{tcolorbox}[boxrule=0pt, frame empty, colback=lightgray, breakable]
\begin{center}
\textbf{#1}
\end{center}
#2
\end{tcolorbox}
}{%
#2}}

\newcommand{\feasibleSet}         {\mathcal{C}}
\newcommand{\objective}           {\mathcal{J}}
\newcommand{\varGeneric}          {\mu}
\newcommand{\varOther}            {\nu}
\newcommand{\varFixed}            {\bar\mu}
\newcommand{\varReference}        {\hat\mu}
\newcommand{\varOptimal}          {\mu^\ast}
\newcommand{\varTangent}          {\xi}
\newcommand{\varTangentReals}     {\upsilon}

\newcommand{\varTangentFixed}     {\bar\varTangent}
\newcommand{\varTangentCoupling}  {\alpha}
\newcommand{\varPlanGeneric}      {\gamma}

\newcommand{\localZero}[1]{\mathbf{0}_{#1}}
\newcommand{\localEquiv}[1]               {=}
\newcommand{\localSubseteq}[1]               {\subseteq} 

\newcommand{\localSum}[1]{\,+_{#1}\,}

\newcommand{\onotation}[1]                   {o\left(#1\right)}

\newif\ifmargincomments
\margincommentsfalse
\ifmargincomments
	
\else
	
\fi

\input{utils/tikzart}
\usepackage[acronym]{glossaries}
\newacronym{acr:amod}{AMoD}{autonomous mobility-on-demand}
\newacronym{acr:av}{AV}{autonomous vehicle}
\newacronym{acr:are}{ARE}{algebraic riccati equation}
\newacronym{acr:dare}{DARE}{discrete-time algebraic riccati equation}
\newacronym{acr:dro}{DRO}{distributionally robust optimization}
\newacronym{acr:dpa}{DPA}{dynamic programming algorithm}
\newacronym{acr:isc}{isc}{inner semicontinuous}
\newacronym{acr:kl}{KL}{ullback--Leibler}
\newacronym{acr:lsc}{lsc}{lower semicontinuous}
\newacronym{acr:osc}{osc}{outer semicontinuous}
\newacronym{acr:pmp}{PMP}{Pontryagin maximum principle}
\newacronym{acr:rl}{RL}{reinforcement learning}
\newacronym{acr:snc}{SNC}{sequentially normally compact}
\newacronym{acr:usc}{usc}{upper semicontinuous}

\boolfalse{proofInMain}
\booltrue{proofInAppendix}
\boolfalse{duplicateStatementInAppendix}
\boolfalse{showLinks}
\boolfalse{showBoxes}

\def\argmax{\mathop{\rm arg\,max}}
\def\argmin{\mathop{\rm arg\,min}}

\theoremstyle{remark}
\newtheorem*{remark*}{Remark}

\numberwithin{equation}{section}

\crefname{figure}{Figure}{Figures}
\Crefname{figure}{Figure}{Figures}

\title{Variational Analysis in the Wasserstein Space}
\author{%
Nicolas Lanzetti$^*$\\
Caltech\\
\texttt{lnicolas@caltech.edu} \\
\And
Antonio Terpin$^*$\\
ETH Z\"urich\\
\texttt{aterpin@ethz.ch}
\And
Florian D\"orfler\\
ETH Z\"urich\\
\texttt{dorfler@ethz.ch}
}

\begin{document}

\maketitle

\footnotefirstpage{$^*$ Equal contribution.}
\footnotefirstpage{%
This research was carried out while all authors were affiliated with ETH Z\"urich.
It was supported by the Swiss National Science Foundation under the NCCR Automation, grant agreement 51NF40\_180545.
}

\begin{abstract}
    \noindent
We study optimization problems whereby the optimization variable is a probability measure. Since the space of probability measures is not a vector space, many classical methods for optimization (e.g., gradients) do not directly apply. Thus, one typically resorts to the abstract machinery of infinite-dimensional analysis or other ad-hoc methodologies, not tailored to the space of probability measures, which however involve projections or rely on convexity-type assumptions. We believe instead that these problems call for a comprehensive methodological framework for calculus in the space of probability measures. 
In this work, we combine ideas from optimal transport, variational analysis, and Wasserstein gradient flows to equip the Wasserstein space (i.e., the space of probability measures endowed with the Wasserstein distance) with a variational structure, both by combining and extending existing results and introducing novel tools. Our theoretical analysis culminates in general necessary optimality conditions.
These conditions (i) resemble the optimality conditions in Euclidean spaces, such as the KKT conditions, (ii) are intuitive, informative, and easy to study, and (iii) yield closed-form solutions or can be used to design computationally attractive algorithms. We accompany our theoretical results with numerous examples and present applications to machine learning, drug discovery, and distributionally robust optimization.
\end{abstract}

\section{Main result}\label{sec:introduction}
This work considers optimization problems of the form
\begin{equation}\label{equation:problem:opt:general}
    \inf_{\varGeneric \in \feasibleSet}\, \objective(\varGeneric),
\end{equation}
where $\feasibleSet \subseteq \Pp{}{\reals^d}$ is a set of admissible \WikipediaLink{https://en.wikipedia.org/wiki/Probability_measure}{probability measures} and $\objective: \Pp{}{\reals^d} \to \realsBar$ is a functional to minimize. This abstract problem setting stems from the observation that numerous fields, including machine learning, robust optimization, and biology, tackle their own version of~\eqref{equation:problem:opt:general}, but with ad-hoc methods that often rely on structural assumptions that do not hold in broader settings.
Despite the recent efforts in the literature~\cite{Lanzetti2022}, these problems still demand a comprehensive theory for the optimization problem~\eqref{equation:problem:opt:general}, which is the subject of this work.
Specifically,
we derive novel \WikipediaLink{https://en.wikipedia.org/wiki/Necessity_and_sufficiency}{necessary first-order optimality conditions} for \eqref{equation:problem:opt:general}, for arbitrary functionals and constraints. These formally resemble the rationales of Euclidean spaces (e.g., \WikipediaLink{https://en.wikipedia.org/wiki/Karush--Kuhn--Tucker_conditions}{Karush--Kuhn--Tucker} and \WikipediaLink{https://en.wikipedia.org/wiki/Lagrange_multiplier}{Lagrange conditions}) and are intuitive, informative, and easy to study.
As a byproduct of our analysis, we translate tools from \WikipediaLink{https://en.wikipedia.org/wiki/Variational_analysis}{variational analysis} (e.g., generalized \WikipediaLink{https://en.wikipedia.org/wiki/Subderivative}{subgradients}, \WikipediaLink{https://en.wikipedia.org/wiki/Normal_cone}{normal cones}, \WikipediaLink{https://en.wikipedia.org/wiki/Tangent_cone}{tangent cones}, etc.) to the Wasserstein space (i.e., the space of probability measures endowed with the \WikipediaLink{https://en.wikipedia.org/wiki/Wasserstein_metric}{Wasserstein distance}).
After illustrating these novel tools in numerous pedagogical examples, we tackle open problems arising in machine learning, drug discovery, and \gls*{acr:dro}, showcasing how our conditions result either in closed-form solutions of \eqref{equation:problem:opt:general} or computationally attractive algorithms.

Our main results are general first-order optimality conditions of~\eqref{equation:problem:opt:general}:

{
\begin{informaltheorem}[First-order optimality conditions%
]
\label{informaltheorem}
If $\varOptimal \in \Pp{}{\reals^d}$ is an optimal solution of \eqref{equation:problem:opt:general} with finite second moment and if a constraint qualification holds, then the ``Wasserstein subgradients'' are ``aligned'' with the constraints at ``optimality'', i.e.,
\[
\localZero{\varOptimal} \in \subgradient{\objective}(\varOptimal) + \normalCone{\feasibleSet}{\varOptimal},
\]
where $\subgradient{\objective}(\varOptimal)$ is the ``Wasserstein subgradient'' of $\objective$ at $\varOptimal$, $\normalCone{\feasibleSet}{\varOptimal}$ is the ``Wasserstein normal cone'' of $\feasibleSet$ at $\varOptimal$ and $\localZero{\varOptimal}$ is a ``null Wasserstein tangent vector'' at $\varOptimal$.
\end{informaltheorem}
}

\begin{figure}[t]
\centering
\begin{minipage}{.45\textwidth}
\centering
\includegraphics[width=.75\linewidth]{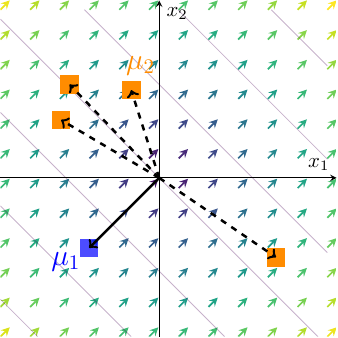}
\end{minipage}
\hfill
\begin{minipage}{.45\textwidth}
\centering
\includegraphics[width=.75\linewidth]{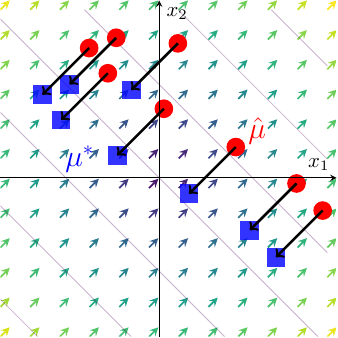}
\end{minipage}
\caption{``Gradients'' are ``aligned'' with the constraints at optimality.
The figure on the left depicts two empirical candidate solutions
$\textcolor{blue}{\varGeneric_1}$ and $\textcolor{orangeDark}{\varGeneric_2}$ for \eqref{equation:problem:opt:general} with $\feasibleSet \coloneqq \{\varGeneric\in\Pp{2}{\reals^2}\st\expectedValue{(x_1, x_2)\sim\varGeneric}{x_1^2+x_2^2}\leq \varepsilon^2\}$ (i.e., bounded second moment) and $\objective(\varGeneric) \coloneqq \expectedValue{(x_1, x_2) \sim \textcolor{blue}{\varGeneric}}{x_1+x_2}$, of which we show the contours and the gradient vector field.
The solid (resp. dashed) black arrows represent the gradient of the constraint function $\expectedValue{(x_1, x_2)\sim\varGeneric}{x_1^2+x_2^2}- \varepsilon^2$ at $\textcolor{blue}{\varGeneric_1}$ (resp. $\textcolor{orangeDark}{\varGeneric_2}$).
Here, $\textcolor{blue}{\varGeneric_1}$ is indeed a candidate optimal solution: The gradient of the objective is aligned with the gradient of the constraint.
For $\textcolor{orangeDark}{\varGeneric_2}$, instead, these two are not aligned. Thus, $\textcolor{orangeDark}{\varGeneric_2}$ cannot be optimal.
The figure on the right shows that $\textcolor{blue}{\varOptimal}$ satisfies the optimality condition for \eqref{equation:problem:opt:general} with $\feasibleSet \coloneqq\{\varGeneric\in\Pp{2}{\reals^2}\st\wassersteinDistance{\varGeneric}{\textcolor{red}{\varReference}}{2}\leq\varepsilon\}$ (i.e., Wasserstein ball centered at $\textcolor{red}{\varReference}$ of radius $\varepsilon$) and $\objective(\varGeneric) \coloneqq \expectedValue{(x_1, x_2) \sim \varGeneric}{x_1 + x_2}$, of which the contours and the gradient vector field are shown. The black arrows connecting $\textcolor{red}{\varReference}$ and $\textcolor{blue}{\varOptimal}$ represent the gradient of the constraint function $\wassersteinDistance{\textcolor{blue}{\varGeneric}}{\textcolor{red}{\varReference}}{2}^2-\varepsilon^2$.
Since $\textcolor{blue}{\varOptimal}$ is optimal, the gradient of the objective and the constraint are aligned at all the ``particles'' of $\textcolor{blue}{\varOptimal}$.
}
\label{fig:first-order-opt-conditions}
\end{figure}

As corollaries of our theorem, we obtain the ``Wasserstein counterparts'' of \WikipediaLink{https://en.wikipedia.org/wiki/Fermat\%27s_theorem_(stationary_points)}{Fermat's rule} in the unconstrained setting (i.e., the gradient vanishes at optimality) and the Lagrange conditions for (in)equality-constrained settings (i.e., the ``gradients'' of the objective and the constraint are ``aligned'' at optimality, see~\cref{fig:first-order-opt-conditions}).

Before diving into variational analysis in the Wasserstein space, we illustrate our optimality conditions by informally studying a simple and accessible version of \eqref{equation:problem:opt:general}.
For $\theta \neq 0$ and $\varepsilon > 0$, consider the problem
\begin{equation}\label{equation:problem:intro:simple}
    \inf_{\varGeneric \in \Pp{2}{\reals^d}}\; \expectedValue{x \sim \varGeneric}{\innerProduct{\theta}{x}}
    \qquad\text{subject to}\qquad
    \expectedValue{x \sim \varGeneric}{\norm{x}^2} \leq \varepsilon^2,
\end{equation}
depicted in \cref{fig:first-order-opt-conditions} for $\theta = \begin{bsmallmatrix}
    1 & 1
\end{bsmallmatrix}^\top$.
To get some intuition, let us restrict to \WikipediaLink{https://en.wikipedia.org/wiki/Dirac_measure}{Dirac deltas} of the form $\delta_x$ for $x\in\reals^d$. Accordingly, \eqref{equation:problem:intro:simple} reduces to
$
    \inf_{\norm{x}^2 \leq \varepsilon^2} \innerProduct{\theta}{x}.
$
This optimization problem can be studied through standard first-order optimality conditions in \WikipediaLink{https://en.wikipedia.org/wiki/Euclidean_space}{Euclidean spaces}.
Since the gradient of the objective $\gradient{}{\innerProduct{\theta}{x}} = \theta$ never vanishes, the optimal solution (if it exists) lies at the \WikipediaLink{https://en.wikipedia.org/wiki/Boundary_(topology)}{boundary}.
We thus seek the \WikipediaLink{https://en.wikipedia.org/wiki/Lagrange_multiplier}{Lagrange multiplier} $\lambda > 0$ such that
\begin{equation}
\label{equation:optimality-conditions:reals}
    0 = \gradient{}{\innerProduct{\theta}{x}} + \lambda\gradient{}{\norm{x}^2} = \theta + 2\lambda x
    \quad\text{and}\quad
    \varepsilon^2 = \norm{x}^2,
\end{equation}
which yields $x = -\varepsilon\frac{\theta}{\norm{\theta}}$ and  $\lambda = \frac{\norm{\theta}}{2\varepsilon}$.
Now back to \eqref{equation:problem:intro:simple}: After basic algebraic manipulations, our main result (stated informally above) tells us that any solution $\varOptimal \in \Pp{2}{\reals^d}$ of \eqref{equation:problem:intro:simple} satisfies the Lagrange-like condition
\begin{equation}
\label{equation:optimality-conditions:ours}
0 = 2\lambda x + \theta\quad\almostEverywhere{\varOptimal}
\quad\text{and}\quad
\varepsilon^2 = \expectedValue{x \sim \varOptimal}{\norm{x}^2} = \frac{\norm{\theta}^2}{4\lambda^2},
\end{equation}
for some $\lambda \geq 0$ constant across the \WikipediaLink{https://en.wikipedia.org/wiki/Support_(measure_theory)}{support} of $\varOptimal$; cf. \cref{fig:first-order-opt-conditions} for $\theta = \begin{bsmallmatrix}
    1 & 1
\end{bsmallmatrix}^\top$.
We conclude that the mass of any candidate solution is necessarily concentrated at $-\varepsilon\frac{\theta}{\norm{\theta}}$.
In particular, our optimality condition in \eqref{equation:optimality-conditions:ours} mirrors its counterpart on $\reals^d$ in \eqref{equation:optimality-conditions:reals}.

\subsection{Introduction to the broader context}

Optimization problems over the space of probability measures are ubiquitous across a variety of fields.

\paragraph{\gls*{acr:dro}.} \gls*{acr:dro} has emerged as a major paradigm for decision-making under uncertainty \cite{rahimian2019distributionally,kuhn2025distributionally}. In \gls*{acr:dro}, one accepts the inherent ambiguity in real-world data and other a priori assumptions in \WikipediaLink{https://en.wikipedia.org/wiki/Uncertainty_quantification}{uncertainty modeling}, where the true underlying probability measure is uncertain and difficult to ascertain, and seeks decisions that are robust against a range of possible probability measures. Thus, \gls*{acr:dro} falls within the scope of \eqref{equation:problem:opt:general} where $\objective$ is a \WikipediaLink{https://en.wikipedia.org/wiki/Risk_measure}{risk measure} \cite{artzner1999coherent,embrechts2022robustness,krokhmal2011modeling} and $\feasibleSet$ is a so-called ambiguity set of probability measures, often defined in terms of the \WikipediaLink{https://en.wikipedia.org/wiki/Kullback–Leibler_divergence}{Kullback-Leibler divergence} \cite{duchi2021statistics,namkoong2017variance,wang2016likelihood} or an \WikipediaLink{https://en.wikipedia.org/wiki/Transportation_theory_(mathematics)}{optimal transport discrepancy} \cite{blanchet2019quantifying,blanchet2022optimal,gao2023distributionally,kuhn2019wasserstein}.

\paragraph{Inverse problems.}
In \WikipediaLink{https://en.wikipedia.org/wiki/Inverse_problem}{inverse problems} one seeks the state of a system given some noisy observations. For instance, Bayesian inference amounts to solving
\begin{equation}\label{equation:bayes_rule_three}
    \inf_{\varGeneric \in \Pp{}{\reals^d}} \expectedValue{x \sim \varGeneric}{\sum_{i = 1}^N-\log(p(z_i|x))} + \kl{\varGeneric}{\varReference}{},
\end{equation}
where $z_i$ for $i \in \{1, \ldots, N\}$ are the observations, $p(z_i|x)$ is the probability of the observation $z_i$ given the state $x \sim \varGeneric$, and $\kullbackLeibler$ is the Kullback-Leibler divergence between a candidate posterior $\varGeneric$ and the prior $\varReference \in \Pp{}{\reals^d}$~\cite{knoblauch2022optimization}. Many \WikipediaLink{https://en.wikipedia.org/wiki/Statistical_inference}{inference problems} (e.g., see \cite[Table 1]{knoblauch2022optimization} and \cite{chu2019probability,rigollet2018entropic}) result from variants of~\eqref{equation:bayes_rule_three}.

\paragraph{Reinforcement learning.}
The dual formulation \cite{puterman1994markov} of the \WikipediaLink{https://en.wikipedia.org/wiki/Reinforcement_learning}{\gls*{acr:rl}} problem seeks the optimal stationary state-action distribution $\varOptimal \in \Delta(\mathcal{S}\times\mathcal{A}) \subseteq \Pp{}{\mathcal{S}\times\mathcal{A}}$, where $\mathcal{S}$ is the state space, $\mathcal{A}$ the action space, and $\Delta(\mathcal{S}\times\mathcal{A})$ the set of stationary state-action distributions compatible with the dynamics (and constraints), that maximizes the expected reward $r: \mathcal{S}\times\mathcal{A} \to \realsBar$,
\begin{equation}
\label{equation:rl}
\varOptimal \in \argmax_{\varGeneric \in \Delta(\mathcal{S}\times\mathcal{A})} \expectedValue{(s, a) \sim \varGeneric}{r(s, a)}.
\end{equation}
Variations of \eqref{equation:rl} yield different problem settings~\cite{greenberg2022efficient,terpin2022trust,wachi2020safe}.

\paragraph{Many others.} A growing number of fields are tackling optimization problems formulated in the space of probability measures, including weather forecasting \cite{fisher2009data}, single-cell perturbation responses \cite{bunne2023single-cell-perturbation,terpin2024learning}, neural network training~\cite{chizat2018global}, generative models \cite{balcerak2025energy}, control of dynamical systems \cite{pilipovsky2024distributionally,terpin2023dynamic}, mean-field control~\cite{bongini2017mean,bonnet2019pontryagin}, and finance \cite{markowitz1952portfolio,rahimian2019distributionally}, among others.

\subsection{Related work}
Our work studies~\eqref{equation:problem:opt:general} through the lens of optimal transport.
The theory of optimal transport, dating back to the seminal work of \WikipediaLink{https://en.wikipedia.org/wiki/Gaspard_Monge}{Monge} \cite{monge1781memoire} and \WikipediaLink{https://en.wikipedia.org/wiki/Leonid_Kantorovich}{Kantorovich} \cite{kantorovich2006masses}, defines a metric, the Wasserstein metric, on the space of probability measures.
While the space of probability measures is not a vector space, which makes most optimization tools in \WikipediaLink{https://en.wikipedia.org/wiki/Banach_space}{Banach spaces} (e.g., \cite{kurcyusz1976existence,luenberger1997optimization}) inapplicable, the theory of optimal transport enables a notion of differentiability, called Wasserstein differentiability, specifically tailored to the Wasserstein space~\cite{ambrosio2005gradient,jko1998,santambrogio2017euclidean}.
Wasserstein differentiability was exploited in~\cite{Lanzetti2022} to derive first-order optimality conditions for \eqref{equation:problem:opt:general} for problems with differentiable objectives and feasible sets given by \WikipediaLink{https://en.wikipedia.org/wiki/Level_set}{sublevel sets} of differentiable functionals.
In the context of optimal control, \cite{bonnet2019pontryagin2,bonnet2021necessary,bonnet2019pontryagin} use Wasserstein differentiability to derive optimality conditions for optimal control problems, whereas \cite{terpin2023dynamic} explores the properties of the dynamic programming algorithm in the space of probability measures in discrete time. The theory has also been applied to derive algorithms to compute Wasserstein barycenters \cite{agueh2011barycenters,chewi2020gradient,panaretos2020invitation}, to analyze over-parametrized neural networks \cite{bach2021gradient,chizat2018global}, approximate inference \cite{frogner2020approx}, and reinforcement learning \cite{zhang2018policy,terpin2022trust}.

Unfortunately, many functionals over the Wasserstein space (including the squared Wasserstein distance itself) fail to be differentiable. This, together with the rigidity of the feasible sets $\feasibleSet$ considered so far, effectively hinders the deployment of these tools in many practical instances.
From a technical standpoint, these limitations are intrinsic in the choice of the variations in the Wasserstein space: Variations induced by the \WikipediaLink{https://en.wikipedia.org/wiki/Pushforward_measure}{pushforward} of (sufficiently regular) \emph{transport maps}, as in~\cite{Lanzetti2022}, are not expressive enough. For instance, the pushforward of an \WikipediaLink{https://en.wikipedia.org/wiki/Empirical_measure}{empirical probability measure} is always empirical (in particular, mass cannot be split).
Thus, while attractive (e.g., perturbations can be captured by well-behaved functions, which form a vector space), a comprehensive theory of calculus requires a more general way of perturbing probability measures.

Intuitively, we can approximate a Dirac delta by \WikipediaLink{https://en.wikipedia.org/wiki/Gaussian_measure}{Gaussians} with vanishing variance. However, there is no transport map describing a variation from the Dirac delta to the approximating Gaussians. In this work, we therefore adopt a different approach and consider perturbations induced by \emph{transport plans}. This more general approach requires us to dive into the generalized notion of Wasserstein subgradient proposed in~\cite[Section 10.3]{ambrosio2005gradient} and, importantly and contrary to the literature, in the generalized tangent space first introduced in~\cite[Section 12]{ambrosio2005gradient}.
Although these perturbations entail significant challenges (e.g., transport plans do not form a vector space), we show in this work that, if judiciously combined with traditional ideas from variational analysis~\cite{rockafellar2009variational,mordukhovich2006variationalI,mordukhovich2006variationalII,mordukhovich2018variational}, they result in general necessary optimality conditions for~\eqref{equation:problem:opt:general}.

Our approach offers several advantages over alternative methods for optimization in the space of probability measures.
First, our analysis is specifically tailored to the Wasserstein space and does not require the introduction of nonnegativity and normalization. For instance, if~\eqref{equation:problem:opt:general} is unconstrained, the corresponding optimality condition simply states that the gradient vanishes at optimality---just like in Euclidean settings.
Second, our optimality conditions hold in full generality, and, in particular, do not rely on convexity-type or linearity assumptions, which in some cases might allow one to study~\eqref{equation:problem:opt:general} through infinite-dimensional linear programming or convex analysis.
Third, analogously to the traditional Karush--Kuhn--Tucker conditions, we demonstrate that our optimality conditions can be used to both solve~\eqref{equation:problem:opt:general} in closed form, and devise efficient numerical methods when a closed-form solution is not available.

\subsection{More details on our contributions}

More specifically, our main contribution consists of four key aspects.
First, we provide a solid foundation for a theory of variational analysis in the Wasserstein space, introducing various concepts such as the generalized subgradient, normal cone, and tangent cone in the Wasserstein space.
For this, we need to resort to very general perturbations induced by transport plans. While invisible to the end user, these perturbations entail working with a tangent space which is not a linear space and carefully dealing with compactness issues.
Second, we provide closed-form and easy-to-use expressions for the subgradients of many functionals and for the normal cone to feasible sets of practical interest.
In particular, we show that the squared Wasserstein distance is not regularly subdifferentiable and only admits a generalized subgradient. This result, which we believe to be of independent interest, also confirms that a general theory of optimality conditions is required already for simple functions (in particular, the distance itself), and not only to cover all corner cases. The key technique to characterize the general subgradient of the squared Wasserstein distance involves its differentiability at regular measures together with approximation arguments, which is a promising technique to explore in future work addressing the differentiability of other functionals not covered in this work.
Third, we derive general first-order optimality conditions for~\eqref{equation:problem:opt:general}. Inspired by classical variational analysis in Euclidean spaces, we prove our result by reformulating~\eqref{equation:problem:opt:general} as an optimization problem over the epigraph of $\objective$. This way, we can establish our optimality conditions under very weak assumptions on the functionals and constraints.
Notably, as we demonstrate with several pedagogical examples, the technical complexity behind our optimality conditions is hidden from the end user, and their deployment is effectively analogous to what one would do in Euclidean spaces.
Fourth, we deploy our results to study a wide variety of optimization problems of the form~\eqref{equation:problem:opt:general} arising in machine learning, drug discovery, and~\gls{acr:dro}. Across all these settings, we show that our optimality conditions both enable novel insights and, when the problem of interest does not admit a closed-form solution, can be used to design effective computational methods.
We believe our tools enable the development of novel algorithms and results in machine learning, robust optimization, and biology, among others.

\section{Subgradients and variational geometry}
\label{section:wasserstein_space}
After recalling preliminaries in measure theory and optimal transport in
\cref{section:ot_wasserstein:background,section:ot_wasserstein:optimal_transport_background}, we present variations in the
Wasserstein space in \cref{subsec:variations}, introduce Wasserstein
subgradients in \cref{subsec:wasserstein_subgradients}, and study the
variational geometry of the Wasserstein space in
\cref{section:optimization:var_geom}.

\subsection{Preliminaries}
\label{section:ot_wasserstein:background}
All the maps considered in this work are tacitly assumed to be Borel, i.e., \WikipediaLink{https://en.wikipedia.org/wiki/Measurable_function}{measurable} with respect to the \WikipediaLink{https://en.wikipedia.org/wiki/Borel_set}{Borel $\sigma$-algebra}.
The set of \WikipediaLink{https://en.wikipedia.org/wiki/Borel_measure}{Borel probability measures} on $\reals^d$ is $\Pp{}{\reals^d}$, and we denote by $\Pp{2}{\reals^d}$ the set of probability distributions with finite \WikipediaLink{https://en.wikipedia.org/wiki/Moment_(mathematics)}{second moment}. We write $\varGeneric \ll \varFixed$ to indicate that $\varGeneric$ is \WikipediaLink{https://en.wikipedia.org/wiki/Absolute_continuity}{absolutely continuous} with respect to $\varFixed$, and we denote by $\Ppabs{2}{\reals^d}$ the set of absolutely continuous probability measures with respect to the Lebesgue measure with finite second moment. The \WikipediaLink{https://en.wikipedia.org/wiki/Support_(measure_theory)}{support} $\supp(\varGeneric) \subseteq \reals^d$ of a probability measure $\varGeneric \in \Pp{}{\reals^d}$ is the closed set
$
\left\{
x \in \reals^d
\st
\varGeneric(U) > 0
\text{ for each neighborhood $U$ of $x$}
\right\}.
$ The identity map on $\reals^d$ is $\identity{\reals^{d}}(x) = x$ and, when clear from the context, we simply write $\identity{}$. The \WikipediaLink{https://en.wikipedia.org/wiki/Gradient}{gradient} of a function $h: \reals^d \to \reals$ at $x \in \reals^d$ is $\gradient{}{h}(x)$, and the \WikipediaLink{https://en.wikipedia.org/wiki/Partial_derivative}{partial derivatives} of a function $c: \reals^d\times\reals^d\to\reals$ at $(x, y) \in \reals^d\times\reals^d$ are denoted $\gradient{x}{c}(x, y)$ and $\gradient{y}{c}(x, y)$, respectively. For two expressions $f(x)$ and $g(x)$, we \WikipediaLink{https://en.wikipedia.org/wiki/Big_O_notation\#Little-o_notation}{write} $f(x) = \onotation{g(x)}$ if, as $g(x) \to 0$, we have $\frac{f(x)}{g(x)} \to 0$.
The \emph{pushforward} of a probability measure $\varGeneric \in \Pp{}{\reals^d}$ through a map $T: \reals^d \to \reals^n$ (cf. \cite[Definition 1.2.2]{figalli2021invitation}), denoted by $\pushforward{T}{\varGeneric}\in\Pp{}{\reals^n}$, is defined by $(\pushforward{T}{\varGeneric})(B) = \varGeneric(T^{-1}(B))$ for all Borel sets $B \subseteq \reals^n$, and it is a probability measure; see \cite[Lemma 1.2.3]{figalli2021invitation}. Then, for any $\pushforward{T}{\varGeneric}$-integrable $\phi: \reals^n \to \reals$, it holds
$\int_{\reals^n} \phi(y) \d(\pushforward{T}{\varGeneric})(y) = \int_{\reals^{d}} \phi(T(x)) \d\varGeneric(x)$; see \cite[Corollary 1.2.6]{figalli2021invitation}.
Finally (cf. \cite[Lemma 1.2.7]{figalli2021invitation}), for any {$T: \reals^d \to \reals^n$ and $S: \reals^n \to \reals^t$ measurable,
$\pushforward{(\comp{T}{S})}{\varGeneric} = \pushforward{S}{(\pushforward{T}{\varGeneric})}$}.

\subsection{Optimal transport}
\label{section:ot_wasserstein:optimal_transport_background}
Given $\varGeneric_1 \in \Pp{2}{\reals^d}$ and $\varGeneric_2 \in \Pp{2}{\reals^n}$, their \WikipediaLink{https://en.wikipedia.org/wiki/Product_measure}{product measure} is $\varGeneric_1\times\varGeneric_2$. We say that $T: \reals^d \to \reals^n$ is a \emph{transport map} from $\varGeneric_1$ to $\varGeneric_2$ if $\pushforward{T}{\varGeneric_1} = \varGeneric_2$ or, equivalently,  $\int_{\reals^n} \phi(y) \d \varGeneric_2(y) = \int_{\reals^n} \phi(y) \d (\pushforward{T}{\varGeneric_1})(y)=\int_{\reals^d}\phi(T(x))\d\varGeneric_1(x)$ for all $\phi \in \Cb{\reals^n}$~\cite[Lemma 1.2.5]{figalli2021invitation}, where $\Cb{\reals^n}$ denotes the space of real-valued \WikipediaLink{https://en.wikipedia.org/wiki/Continuous_function}{continuous} \WikipediaLink{https://en.wikipedia.org/wiki/Bounded_function}{bounded} functions on $\reals^n$. For some $i \in \naturals$, $1 \leq i \leq n$, we denote by $\projection{i}{}: (\reals^{m})^n \to \reals^{m}$ the projection map on the $i^\mathrm{th}$ component, i.e., $\projection{i}{}(x_1, x_2, \ldots, x_n) = x_i$. That is, $\projection{i}{}$ always denotes the projection onto the $i^\mathrm{th}$ factor of the relevant product space; its precise domain and codomain are clear from the context. We use the projection map to obtain the \WikipediaLink{https://en.wikipedia.org/wiki/Marginal_distribution}{\emph{marginals}} of a probability measure via pushforward, $\pushforward{\projection{i}{}}{\varGeneric_1}$. A \emph{transport plan} between $\varGeneric_1$ and $\varGeneric_2$ is a probability measure $\varPlanGeneric \in \Pp{2}{\reals^d \times \reals^n}$ so that $\pushforward{\projection{1}{}}{\varPlanGeneric} = \varGeneric_1$ and $\pushforward{\projection{2}{}}{\varPlanGeneric} = \varGeneric_2$. A transport map may not exist between $\varGeneric_1$ and $\varGeneric_2$ (e.g., when $\varGeneric_1 = \delta_{x}$ and $\varGeneric_2$ is a Gaussian distribution), but a transport plan always does (e.g., the product measure $\varPlanGeneric = \varGeneric_1 \times \varGeneric_2$). We collect them in the set of \emph{transport plans} $\setPlans{\varGeneric_1}{\varGeneric_2}$. Given a lower semicontinuous function $c: \reals^d\times\reals^n\to\realsBar \coloneqq \reals\cup\{\pm\infty\}$, the optimal transport problem reads:
\begin{equation}
\label{definition:optimaltransport}
\wassersteinDistance{\varGeneric_1}{\varGeneric_2}{c}
\coloneqq
\min_{\varPlanGeneric \in \setPlans{\varGeneric_1}{\varGeneric_2}}
\int_{\reals^d\times\reals^n}
c(x, y)
\d\varPlanGeneric(x, y).
\end{equation}
The celebrated (quadratic) Wasserstein distance is a special case of \eqref{definition:optimaltransport} (cf. \cite[Definition 3.1.3]{figalli2021invitation}), with $n = d$ and $c(x, y) = \norm{x - y}^2$:
\emphBox{Type $2$ Wasserstein distance}{
\begin{equation}
\label{definition:wasserstein-distance}
\wassersteinDistance{\varGeneric_1}{\varGeneric_2}{2}
\coloneqq
\left(
\min_{\varPlanGeneric \in \setPlans{\varGeneric_1}{\varGeneric_2}}
\int_{\reals^d\times\reals^d}
\norm{x - y}^2
\d\varPlanGeneric(x, y)
\right)^{\frac{1}{2}},
\end{equation}
}
where $\norm{\cdot}$ is the standard \WikipediaLink{https://en.wikipedia.org/wiki/Norm_(mathematics)}{Euclidean norm} in $\reals^d$.
We write $\setPlansOptimalCost{\varGeneric_1}{\varGeneric_2}{c}$ and $\setPlansOptimal{\varGeneric_1}{\varGeneric_2}$ for the set of minimizers of \eqref{definition:optimaltransport} and \eqref{definition:wasserstein-distance}, respectively.
The Wasserstein distance is a distance on $\Pp{2}{\reals^d}$~\cite[Section 6]{villani2009optimal}.
We define the Wasserstein ball of radius $\varepsilon$ as
$
\wassersteinBall{\varFixed}{\varepsilon}{2}
\coloneqq
\left\{
\varGeneric \in \Pp{2}{\reals^d}
\st
\wassersteinDistance{\varGeneric}{\varFixed}{2} < \varepsilon
\right\}.
$

Throughout the work, we use two notions of \WikipediaLink{https://en.wikipedia.org/wiki/Convergence_of_random_variables}{convergence} for probability measures. A \WikipediaLink{https://en.wikipedia.org/wiki/Sequence}{sequence} $(\varGeneric_n)_{n \in \naturals} \subseteq \Pp{}{\reals^d}$ (i) converges narrowly to $\varGeneric \in \Pp{}{\reals^d}$, denoted by $\varGeneric_n \convergenceNarrow \varGeneric$, if for all $\phi \in \Cb{\reals^d}$ we have $\int_{\reals^d} \phi(x) \d\varGeneric_n(x)
\to
\int_{\reals^d} \phi(x) \d\varGeneric(x)$ (cf. \cite[Definition 2.1.5]{figalli2021invitation}) and (ii) {converges in Wasserstein to $\varGeneric \in \Pp{2}{\reals^d}$,
$
\varGeneric_n \convergenceWasserstein \varGeneric
$, if $
\lim_{n \to \infty}
\wassersteinDistance{\varGeneric_n}{\varGeneric}{2}
= 0
$. The narrow \WikipediaLink{https://en.wikipedia.org/wiki/Topology}{topology} is \WikipediaLink{https://en.wikipedia.org/wiki/Comparison_of_topologies}{weaker} than the Wasserstein topology on $\Pp{2}{\reals^d}$. Indeed, by~\cite[Proposition 7.1.5]{ambrosio2005gradient},
$\varGeneric_n \convergenceWasserstein \varGeneric$ if and only if $\varGeneric_n \convergenceNarrow \varGeneric$ and $\int_{\reals^d}\norm{x}^2\d\varGeneric_n(x)\to\int_{\reals^d}\norm{x}^2\d\varGeneric(x)$ or, equivalently, $\int_{\reals^d}\phi(x)\d\varGeneric_n(x)\to\int_{\reals^d}\phi(x)\d\varGeneric(x)$ for all real-valued continuous $\phi$ with $\abs{\phi(x)}\leq C(1+\norm{x}^2)$.}

\subsection{Variations in the Wasserstein space}
\label{subsec:variations}

\begin{figure}
    \centering
    \begin{minipage}{.325\textwidth}
    \centering
    \includegraphics[width=.45\linewidth]{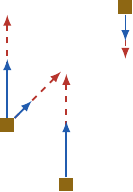}
    \end{minipage}
    \hfill
    \begin{minipage}{.325\textwidth}
    \centering
    \includegraphics[width=.95\linewidth]{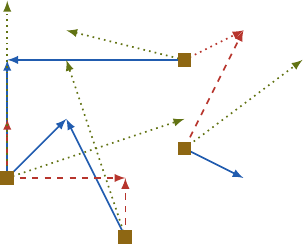}
    \end{minipage}
    \hfill
    \begin{minipage}{.325\textwidth}
    \centering
    \includegraphics[width=.95\linewidth]{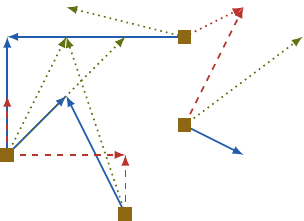}
    \end{minipage}
    \caption{{Wasserstein geometry}.
    On the left, we consider example variations for $\textcolor{ETHbronze}{\varGeneric} \in \Pp{2}{\reals^d}$ (denoted by a \textcolor{ETHbronze}{brown square}), $\textcolor{ETHred}{\varTangent_1}$ in \textcolor{ETHred}{dotted red} and $\textcolor{ETHblue}{\varTangent_2}$ in \textcolor{ETHblue}{solid blue}. They are scaled versions of each other: $\textcolor{ETHred}{\varTangent_1} = 2\textcolor{ETHblue}{\varTangent_2}$. The other two images show two possible sums (in \textcolor{ETHgreen}{dotted green}) of the variations $\textcolor{ETHred}{\varTangent_1}$ (in \textcolor{ETHred}{dotted red}) and $\textcolor{ETHblue}{\varTangent_2}$ (in \textcolor{ETHblue}{solid blue}) at $\textcolor{ETHbronze}{\varGeneric} \in \Pp{2}{\reals^d}$: Each sum results from a different coupling $\varTangentCoupling\in\setPlansCommonMarginal{\textcolor{ETHred}{\varTangent_1}}{\textcolor{ETHblue}{\varTangent_2}}{1}$ between $\textcolor{ETHred}{\varTangent_1}$ and $\textcolor{ETHblue}{\varTangent_2}$ (i.e., a different coupling between the arrows starting from the same square). The coupling resulting in the sum on the left is the one maximizing the inner product (cf. \eqref{definition:geometry:inner_product}) of the two variations (and, equivalently, the one attaining the minimum in the distance between the two variations; cf. \eqref{definition:geometry:distance}).
    }
    \label{fig:variations}
\end{figure}

In the Euclidean space $\reals^d$, a variation at $x\in\reals^d$ can be interpreted as an ``arrow'' $\varTangentReals\in\reals^d$ rooted at $x$. In the same spirit, in the Wasserstein space, a variation at $\varFixed \in \Pp{2}{\reals^d}$ is a ``(weighted) collection of arrows'' for each point in the support of $\varFixed$~\cite[Chapter 12]{ambrosio2005gradient}. Formally, a \emph{variation at $\varFixed \in \Pp{2}{\reals^d}$} is a probability measure $\varTangent\in\Pp{2}{\reals^d\times\reals^d}$ whose first marginal is $\varFixed$ (i.e., $\pushforward{\projection{1}{}}{\varTangent} = \varFixed$) (see~\cref{fig:variations}, left).
We denote the space of all variations at $\varFixed$ by $\PpLocal{2}{\reals^d\times\reals^d}{\varFixed} \coloneqq \{\varTangent\in\Pp{2}{\reals^d\times\reals^d} \st \pushforward{\projection{1}{}}{\varTangent} = \varFixed\}$. We can disintegrate (cf. \cite[Theorem 5.3.1]{ambrosio2005gradient}) $\varTangent$ to obtain a collection $\{\varTangent_x\}_{x\in\reals^d}\subseteq\Pp{2}{\reals^d}$ where each $\varTangent_x$ denotes the probability measure on the variations, i.e., the ``(weighted) collection of arrows'' $\almostEverywhere{\varFixed}$
Our definition of variations is in line with the variational analysis literature \cite{rockafellar2009variational,mordukhovich2006variationalI}, where tangent vectors are expressed as the limits of scaled difference vectors.
The distance between two variations $\varTangent_1, \varTangent_2 \in \PpLocal{2}{\reals^d\times\reals^d}{\varFixed}$ can be defined in terms of the distance between their (weighted) arrows. To do so, we have to account for the fact that different arrows starting from the same point can be coupled in different ways (cf. \cref{fig:variations}, middle and right).
We express this coupling as a measure $\varTangentCoupling\in\Pp{2}{\reals^d\times\reals^d\times\reals^d}$ so that $\pushforward{\projection{12}{}}{\varTangentCoupling} = \varTangent_1$ and $\pushforward{\projection{13}{}}{\varTangentCoupling} = \varTangent_2$, where $(\bar{x},\varTangentReals_1,\varTangentReals_2) \in \supp(\varTangentCoupling)$ if and only if the arrows $\varTangentReals_1, \varTangentReals_2\in\reals^d$ are anchored at $\bar{x}$, and $\varTangentReals_1$ and $\varTangentReals_2$ are coupled. Note that $\pushforward{\projection{1}{}}{\varTangentCoupling} = \pushforward{\projection{1}{}}{\varTangent_1} = \pushforward{\projection{1}{}}{\varTangent_2} = \varFixed$, since $\varTangent_1, \varTangent_2$ are variations at $\varFixed$. We refer to $\varTangentCoupling$ as a \emph{coupling} of the two variations $\varTangent_1$ and $\varTangent_2$, reserving the term \emph{transport plan} for the analogous object between two probability measures: a transport plan couples two measures, whereas a coupling couples two variations sharing the common anchor $\bar x$.
Then, the distance between $\varTangent_1$ and $\varTangent_2$ with common anchor $\varFixed$ amounts to the minimum distance that can be obtained among all such couplings (cf. \cref{fig:variations}, middle):
\begin{align}
    \label{definition:geometry:distance}
    \wassersteinDistance{\varTangent_1}{\varTangent_2}{\varFixed}
    &\coloneqq
    \left(
    \min_{\varTangentCoupling \in \setPlansCommonMarginal{\varTangent_1}{\varTangent_2}{1}}
    \int_{\reals^d\times\reals^d\times\reals^d}
    \norm{\varTangentReals_1-\varTangentReals_2}^2
    \d\varTangentCoupling(\bar x, \varTangentReals_1, \varTangentReals_2)
    \right)^{\frac{1}{2}},
\end{align}
where $\setPlansCommonMarginal{\varTangent_1}{\varTangent_2}{1}$ is the set of all couplings, as defined above, between $\varTangent_1$ and $\varTangent_2$:
    \[\setPlansCommonMarginal{\varTangent_1}{\varTangent_2}{1}
    \coloneqq
    \left\{
    \varTangentCoupling \in \Pp{2}{\reals^d\times\reals^d\times\reals^d}
    \st
    \pushforward{\projection{12}{}}{\varTangentCoupling} = \varTangent_1,
    \pushforward{\projection{13}{}}{\varTangentCoupling} = \varTangent_2
    \right\}.\]
In view of \cite[Proposition 12.4.6]{ambrosio2005gradient}, $\wassersteinDistance{\varTangent_1}{\varTangent_2}{\varFixed} = 0$ implies $\wassersteinDistance{\varTangent_1}{\varTangent_2}{2} = 0$.
Similarly, we can also define the \WikipediaLink{https://en.wikipedia.org/wiki/Dot_product}{inner product} and norm:
\begin{align}
    \label{definition:geometry:inner_product}
    \innerProduct{\varTangent_1}{\varTangent_2}_{\varFixed}
    &\coloneqq
    \max_{\varTangentCoupling \in \setPlansCommonMarginal{\varTangent_1}{\varTangent_2}{1}}
    \int_{\reals^d\times\reals^d\times\reals^d}
    \innerProduct{\varTangentReals_1}{\varTangentReals_2}
    \d\varTangentCoupling(\bar x, \varTangentReals_1, \varTangentReals_2),
    \\
    \label{definition:geometry:norm}
    \norm{\varTangent_i}_{\varFixed}
    &\coloneqq
    \left(
    \int_{\reals^d\times\reals^d}
    \norm{\varTangentReals_i}^2
    \d\varTangent_i(\bar x, \varTangentReals_i)\right)^{\frac{1}{2}}, \quad i \in \{1, 2\}.
\end{align}
With these definitions, the Euclidean-like identity
$
    \wassersteinDistance{\varTangent_1}{\varTangent_2}{\varFixed}^2
    =
    \norm{\varTangent_1}_{\varFixed}^2
    -
    2\innerProduct{\varTangent_1}{\varTangent_2}_{\varFixed}
    +
    \norm{\varTangent_2}_{\varFixed}^2
$, which we exploit in our proofs, holds. The identity reveals the motivation behind the $\min$ and $\max$ in \eqref{definition:geometry:distance} and \eqref{definition:geometry:inner_product}: They are the only two terms involving couplings and they appear with opposite signs. Intuitively, to minimize the distance between two variations we need to ``maximally align'' the arrows at each particle; cf. \cref{fig:variations}.

In Euclidean spaces, given a variation $\varTangentReals$ at a point $x$, we can construct a new point by altering $x$ according to $\varTangentReals$: $y = x + \varTangentReals$. Analogously, given a variation $\varTangent$ at a probability measure $\varFixed$, we obtain a new probability measure by transporting, for each $x \in \supp(\varFixed)$, mass along the variations $\varTangentReals$, weighted according to $\varTangent$: $\varGeneric = \pushforward{(\projection{1}{} + \projection{2}{})}{\varTangent}$. However, unlike the Euclidean counterpart where the variation displacing $x$ to $y$ is $\varTangentReals = y - x$, in the Wasserstein space we have multiple ways of connecting $\varFixed$ and $\varGeneric$, described by the transport plans $\setPlans{\varFixed}{\varGeneric}$. Each of these $\varPlanGeneric \in \setPlans{\varFixed}{\varGeneric}$ induces a variation $\varTangent = \pushforward{(\projection{1}{}, \projection{2}{} - \projection{1}{})}{\varPlanGeneric}$.
{Moreover, the space of variations is not a vector space. However, we can equip it with an ``almost linear'' structure (see \cref{fig:variations} for an intuition):
\includeElement{definition}{def_almost_linear_structure_wass}{true}{false}
Another difference is that while $\varTangentReals$ describes a \WikipediaLink{https://en.wikipedia.org/wiki/Geodesic}{geodesic} in the Euclidean space, i.e., $\norm{\varTangentReals} = \norm{y - x}$ and $\norm{\varepsilon\varTangentReals} = \varepsilon\norm{y - x}$ for $\varepsilon > 0$, this is not the case for the generic variation $\varTangent$: $\norm{\varTangent}_{\varFixed} \geq \wassersteinDistance{\varFixed}{\varGeneric}{2}$ and $\norm{\pushforward{(\projection{1}{}, \varepsilon(\projection{2}{} - \projection{1}{}))}{\varPlanGeneric}}_{\varFixed} \geq \varepsilon\wassersteinDistance{\varFixed}{\varGeneric}{2}$.} In light of this, to define a meaningful \emph{tangent space}, we consider only the variations that for sufficiently small ``scalings'' (i.e., for sufficiently small $\bar\varepsilon$ below) are optimal and, thus, describe geodesics in the Wasserstein space:
The \emph{tangent space} $\tangentCone{\Pp{2}{\reals^d}}{\varFixed}$ of $\Pp{2}{\reals^d}$ at $\varFixed\in\Pp{2}{\reals^d}$ is
the closure with respect to $\wassersteinDistance{}{}{\varFixed}$ of
\begin{equation}
    \label{definition:tangent_space}
    \begin{aligned}
    \tangentSpaceInterior{\Pp{2}{\reals^d}}{\varFixed}
    \coloneqq
    \biggl\{
    \varTangent \in \PpLocal{2}{\reals^d\times\reals^d}{\varFixed}
    \st
    &\exists \bar{\varepsilon} > 0, \forall\varepsilon \in [0, \bar{\varepsilon}),
    \pushforward{(\projection{1}{}, \projection{1}{} + \varepsilon\projection{2}{})}{\varTangent}
    \in
    \setPlansOptimal{\varFixed}{\pushforward{(\projection{1}{} + \varepsilon\projection{2}{})}{\varTangent}}
    \biggr\},
    \end{aligned}
\end{equation}
a definition first introduced in \cite[Chapter 12]{ambrosio2005gradient}.
In~\cref{proposition:prop_sum_scale_well_defined,proposition:prop_wasserstein_geom}, we discuss several properties of this ``almost linear'' structure of the space of variations.

\begin{remark*}
The sum of $n$ elements $\varTangent_1,\ldots,\varTangent_n$ results from iteratively applying \cref{definition:def_almost_linear_structure_wass}\ref{item:local-sum}, and we write $\varTangent_1\localSum{\varTangentCoupling}\ldots\localSum{\varTangentCoupling}\varTangent_n=\pushforward{(\projection{1}{},\projection{2}{}+\ldots+\projection{n+1}{})}{\varTangentCoupling}$ for some \[{\varTangentCoupling\in\setPlansCommonMarginal{\varTangent_1,\ldots}{\varTangent_n}{1}
    \coloneqq
    \left\{
    \varTangentCoupling \in \Pp{2}{\reals^d\times(\reals^d)^n}
    \st
    \pushforward{\projection{1,i+1}{}}{\varTangentCoupling} = \varTangent_i,
    i \in \{1, \ldots, n\}
    \right\}}.\]
\end{remark*}

\paragraph{Comparison with the literature.}

Most of the literature uses the tangent space
\begin{multline}\label{eq:tangent_space:literature}
    \closure{\left\{
    \gradient{}\varphi \st \varphi\in\mathrm{Cyl}({\reals^d})
    \right\}}{\Lp{2}{\reals^d}{\varFixed}}
    =
    {\closure{\left\{
    \varepsilon(T - \identity{}) \st \pushforward{(\identity{}, T)}{\varFixed} \in \setPlansOptimal{\varFixed}{\pushforward{T}{\varFixed}}, \varepsilon > 0
    \right\}}{\Lp{2}{\reals^d}{\varFixed}},}
\end{multline}
where $\mathrm{Cyl}({\reals^d})$ is the set of smooth cylindrical functions, which in finite spaces is homeomorphic to $\Ccinf{\reals^d}$, the space of compactly-supported, smooth functions; cf. \cite[Chapter 5]{ambrosio2005gradient}. This (simplified) tangent space was
first introduced in \cite[Chapter 8]{ambrosio2005gradient} (see, e.g.,~\cite{Lanzetti2022,bonnet2019pontryagin,bonnet2019pontryagin2,bonnet2021necessary} and~\cite[Theorem 8.5.1]{ambrosio2005gradient} for the proof of the equality in~\eqref{eq:tangent_space:literature}) and can be embedded in \eqref{definition:tangent_space} since for $\gradient{}\varphi$ with $\varphi\in\Ccinf{\reals^d}$, the tangent vector $\varTangent=\pushforward{(\identity,\gradient{}\varphi)}\varFixed$ belongs to $\tangentCone{\Pp{2}{\reals^d}}{\varFixed}$ of $\Pp{2}{\reals^d}$. In fact, (i) its first marginal is $\varFixed$ and (ii) $\pushforward{(\projection{1}{}, \projection{1}{} + \varepsilon\projection{2}{})}{\varTangent}=\pushforward{(\identity{},\identity{}+\varepsilon\gradient{}\varphi)}{\varFixed}$ is, for $\varepsilon$ sufficiently small, a transport plan. The latter follows from $\identity{}+\varepsilon\gradient{}\varphi$ being the gradient of a convex function and, thus, an optimal transport map \cite[Theorem 1.48]{santambrogio2015optimal}.
While attractive for its simplicity (it is a vector space), the tangent space~\eqref{eq:tangent_space:literature} limits the perturbations to transport maps. For absolutely continuous measures, this comes with no loss of generality in view of \WikipediaLink{https://en.wikipedia.org/wiki/Polar_factorization_theorem}{Brenier's theorem} \cite{brenier1991polar}. However, for empirical probability measures (such as those arising in data-driven applications, where one has only access to a dataset consisting of a finite number of observations), the restriction to transport maps effectively limits the possible perturbations (e.g., a Dirac delta can only be transported to another Dirac delta by a transport map).
For these reasons, we instead adopt the tangent space introduced in \cite[Chapter 12]{ambrosio2005gradient}.

\subsection{Wasserstein subgradients}\label{subsec:wasserstein_subgradients}
We now define the \emph{regular} and \emph{general} subgradients, along the lines of \cite[Definition 8.3]{rockafellar2009variational} and \cite[Section 10]{ambrosio2005gradient}, for the Wasserstein space. The definition is inspired by the Euclidean setting, where $\bar\varTangentReals$ is a subgradient of $f:\reals^d\to\realsBar$ at $\bar x\in\reals^d$ if $f(x)-f(\bar x)\geq \innerProduct{\bar\varTangentReals}{x-\bar x}+\onotation{\norm{x-\bar x}}$ for all $x\in\reals^d$. That is, if $\bar\varTangentReals$ provides a linear lower bound for $f$ at $\bar x$.
General subgradients are then defined as the limits of regular subgradients.
In the Wasserstein space, we can proceed analogously. In particular, we use the (local) inner product~\eqref{definition:geometry:inner_product} and replace the displacement ``$x-\bar x$'' with $\pushforward{(\pi_1, \pi_2 - \pi_1)}{\varPlanGeneric}$, where $\varPlanGeneric\in\setPlansOptimal{\varFixed}{\varGeneric}$ is an optimal transport plan between $\varFixed$ and $\varGeneric$ (which reduces to $\pushforward{(\identity{},T-\identity{})}{\varFixed}$ when $\varPlanGeneric$ is induced by a map $T$):

\includeElement{definition}{def_sub_sup_differentiability}{true}{false}

We consider convergence in Wasserstein for the general subgradient since the first marginals may differ (we consider sequences $\varGeneric_n\convergenceWasserstein\varFixed$). For this reason, the general subgradient may not be in the tangent space, but it has finite second moment.
The regular subgradient, instead, is in the tangent space by definition.
Similarly, we can then define differentiability:

\includeElement{definition}{def_differentiability}{true}{false}

Intuitively, \cref{definition:def_differentiability} states that $\objective$ is differentiable if there is a regular subgradient that is also a ``supergradient'' in the classical sense; cf. \cite[Definition 2.7]{Lanzetti2022}. Wasserstein subgradients enjoy similar properties to their Euclidean analogues:

\includeElement{proposition}{prop_differentiability_properties}{true}{proofInMain}

In particular, \cref{proposition:prop_differentiability_properties}\ref{item:gradient_uniqueness} allows us to talk about \emph{the} gradient. The (general) subgradient might not be a singleton, even for a differentiable functional; see~\cite[Section 8.B]{rockafellar2009variational} for an example in $\reals^d$.
\cref{proposition:prop_differentiability_properties}\ref{item:gradient_characterization}, instead, is the analogue of $f(y)-f(x)=\gradient{} f(x)^\top (y-x) + \onotation{\norm{x-y}}$: Namely, gradients provide a ``linear approximation'' of $\objective$ from below, with errors that are \WikipediaLink{https://en.wikipedia.org/wiki/Big_O_notation\#Little-o_notation}{small} in terms of the Wasserstein distance. In contrast to \eqref{definition:geometry:inner_product}, in \eqref{item:gradient_characterization:2} the maximization over the couplings coincides with the minimization over the same set up to a $\onotation{\wassersteinDistance{\varFixed}{\varGeneric}{2}}$ term, an important property of differentiable functions.
Finally, \cref{proposition:prop_differentiability_properties}\ref{item:gradient_induced_by_map} shows that the gradient, whenever it exists, is always ``deterministic'' (i.e., induced by a transport map) even at measures that are not absolutely continuous (e.g., the empirical measures arising in data-driven applications). Perhaps surprisingly, this holds although our framework deliberately admits subgradients that are genuine transport plans rather than maps (indeed, as we discuss shortly, those of the squared Wasserstein distance are of this form): differentiability nonetheless forces the gradient onto the graph of a map, recovering ``for free'' the map-induced form imposed by definition in~\cite{gangbo2019differentiability,Lanzetti2022}.

\paragraph{Comparison with the literature.}

Our definition of regular subgradient (\cref{definition:def_sub_sup_differentiability}\ref{item:subgradient:regular}) coincides with the definition of extended \WikipediaLink{https://en.wikipedia.org/wiki/Fréchet_derivative}{Fréchet subdifferential}, introduced in~\cite[Definition 10.3.1]{ambrosio2005gradient} and employed, for instance, in~\cite{bonnet2019pontryagin}.
It is a weaker notion of subdifferentiability compared to the ones in~\cite{gangbo2019differentiability,Lanzetti2022}, where subgradients are of the form $\pushforward{(\identity{},\phi)}\varGeneric{}$ for some Borel map $\phi:\reals^d\to\reals^d$.
From an optimization perspective, these definitions are however not satisfactory since, as we will discuss shortly, the squared Wasserstein distance itself fails to be subdifferentiable at measures that are not absolutely continuous.
Inspired by classical variational analysis in Euclidean spaces~\cite{rockafellar2009variational}, we therefore introduce the definition of general subgradients (\cref{definition:def_sub_sup_differentiability}\ref{item:subgradient}).
In the case of variations (and so subdifferentials) induced by transport maps, a similar definition has appeared in~\cite[Definition 11.1.5]{ambrosio2005gradient}.
Besides the more restrictive variations,~\cite[Definition 11.1.5]{ambrosio2005gradient} also uses different notions of convergence and is formulated only for the so-called strong subdifferential~\cite[Definition 10.3.1]{ambrosio2005gradient}, and is thus insufficient for general optimality conditions.
To the best of our knowledge, our definition is therefore novel.

\subsubsection{Why is this machinery needed?}
We argue that our level of sophistication is required for at least two reasons.
First, already in Euclidean spaces, variational analysis is needed to deal with nondifferentiable and nonconvex settings; the (general) subgradient usually appears in the necessary conditions, while the horizon subgradient (which we shall discuss in \cref{section:variational-geometry:epigraphs}) appears in the constraint qualification (e.g., see~\cite[Theorem 8.15]{rockafellar2009variational} or~\cite[Section 5.1]{mordukhovich2006variationalII}). There is no reason to expect the situation to simplify in the more general space of probability measures.
Second, the squared Wasserstein distance itself fails to be regularly subdifferentiable.
The existing theory of gradient flows \cite{ambrosio2005gradient,Lanzetti2022} bypasses this complexity by resorting to absolutely continuous measures. However, this introduces a substantial practical and theoretical limitation since it excludes, among others, many data-driven settings.
By contrast, we show below that the general subgradient is always nonempty.

\paragraph{The squared Wasserstein distance is not regularly subdifferentiable.}
To start, we study the regular subdifferentiability properties of the squared Wasserstein distance:

\includeElement{proposition}{prop_wassersteindistance_regular_subgradient_must_be_plan}{true}{proofInMain}

When the set of optimal transport plans consists of a unique transport plan induced by an optimal transport map $T:\reals^d\to\reals^d$, then $\gradient{}\objective(\varFixed)=\pushforward{(\identity{},\identity{}-T)}{\varFixed}$, which is formally analogous to the gradient of $x\mapsto\frac{1}{2}\norm{x-\hat x}^2$ being $x-\hat x$ in Euclidean spaces (with ``$\identity{}$ as $x$'' and ``$T$ as $\hat x$''). However, and perhaps surprisingly,
\cref{proposition:prop_wassersteindistance_regular_subgradient_must_be_plan} also indicates that the squared Wasserstein distance from a reference measure $\varReference$ is not regularly subdifferentiable at all measures $\varFixed$ for which there exist multiple optimal transport plans between $\varFixed$ and $\varReference$.
Instead, it is regularly subdifferentiable (and thus differentiable) whenever there is a unique optimal transport plan between $\varFixed$ and $\varReference$ and this transport plan is induced by an optimal transport map---a result first established, together with regular subdifferentiability of $-\objective$, in~\cite[Theorem 10.2.6]{ambrosio2005gradient}.
If $\varFixed$ is absolutely continuous there is a unique optimal transport plan induced by an optimal transport map between $\varFixed$ and $\varReference$; cf. \cite{brenier1991polar}. Thus, the squared Wasserstein distance is differentiable at absolutely continuous measures.
This result is sharp: Even when a unique optimal transport plan exists, the squared Wasserstein distance may fail to be regularly subdifferentiable, as we illustrate in the next example.

\includeElement{example}{ex_non_subdifferentiability_regular_wasserstein_singletonplans}{true}{false}

The absence of regular subgradients prompts us to study the general subgradient of the squared Wasserstein distance. We will do so in the next section.

\paragraph{The (general) subgradient of the squared Wasserstein distance.}
The subdifferentiability of the squared Wasserstein distance at absolutely continuous probability measures (cf. \cref{proposition:prop_wassersteindistance_regular_subgradient_must_be_plan}\ref{item:prop_wassersteindistance_regular_subgradient_must_be_plan:differentiable}) fuels the hope that the general subgradient can be studied via a mollification argument, i.e., by approximating each probability measure via absolutely continuous ones. We illustrate the intuition by approximating a Dirac delta at 0 with Gaussian probability measures with vanishing variance:

\includeElement{example}{ex_subdifferentiability_via_gaussians_wasserstein}{true}{false}

Since $\Ppabs{2}{\reals^d}$ is dense in $\Pp{2}{\reals^d}$ \cite[Theorem 2.2.7]{panaretos2020invitation}, we can always construct a sequence $(\varGeneric_n)_{n \in \naturals} \subseteq \Ppabs{2}{\reals^d}$ with $\varGeneric_n\convergenceWasserstein\varFixed$, compute $\varTangent_n \in \subgradientRegular{\objective}(\varGeneric_n)$, which exists by absolute continuity of $\varGeneric_n$, and study the limit of $\varTangent_n$. One way to construct such $\varGeneric_n$ is via convolution with a Gaussian kernel; see \cite[Lemma 7.1.10]{ambrosio2005gradient}.
Overall, this procedure yields a sufficient characterization (for the sake of necessary optimality conditions) of the (general) subgradient of the squared Wasserstein distance:
\includeElement{proposition}{prop_wassersteindistance_general_subgradient}{true}{proofInMain}

In particular, \cref{proposition:prop_wassersteindistance_general_subgradient} establishes subdifferentiability of the squared Wasserstein distance at \emph{all} probability measures $\varFixed \in \Pp{2}{\reals^d}$.

\paragraph{Why not mollify the optimization problem?}
{General subgradients are defined via a mollification argument: A variation is a general subgradient if there are regular subgradients at nearby points that converge to it in Wasserstein.
As established in~\cite{rockafellar2009variational,mordukhovich2018variational} for Euclidean spaces and in the remainder of the paper for the Wasserstein space, general subgradients lead to a coherent framework for calculus and optimization.
However, one might wonder why not \emph{mollify} the optimization problem itself, thereby removing the nonsmoothness, and then solve the resulting smooth problem in the hope that its solutions converge to those of the original nonsmooth problem.
Unfortunately, as we illustrate in the next example, this approach generally fails even in the Euclidean space and, thus, cannot be expected to work in the Wasserstein space either.
}

\begin{example}[Non-convergence of mollified optimization problems]
\label{example:non-convergence-mollification}Consider the optimization problem
    \begin{center}
    \begin{minipage}{.5\linewidth}
        \begin{equation*}
        \min_{x \in \reals} x
        \quad\text{subject to}\quad
        \sqrt{\abs{x}} - x \leq 0.
    \end{equation*}
    \end{minipage}
    \hfill
    \begin{minipage}{.4\linewidth}
    \centering
    \begin{tikzpicture}
  \begin{axis}[
      width=6cm,height=4cm,
      xmin=-1,xmax=2, ymin=-.25,ymax=2,
      axis lines=middle, ticks=none,
      xlabel={$x$}, ylabel={$y$},
      domain=-2:3, samples=200,
      legend style={draw=none, align=left,fill=none,font=\scriptsize,at={(1.02,0.5)},anchor=west}
    ]
    \addplot[very thick] {x};
    \addplot[very thick,dotted] {sqrt(abs(x))};
    \addplot[blue,ultra thick,dashed] coordinates {(-0.05,0) (0.1,0)};
    \addplot[blue,ultra thick,dashed] coordinates {(1,0) (2.75,0)};
    \legend{$x$,$\sqrt{\abs{x}}$,$\feasibleSet$}

  \end{axis}
\end{tikzpicture}
    \end{minipage}
    \end{center}

{
    The constraint $\sqrt{\abs{x}} - x \leq 0$ implies that only nonnegative $x$ are feasible and, thus, the global minimizer is trivially $x = 0$. The function $f(x) = x$ is differentiable everywhere, but the constraint $0 \geq g(x) = \sqrt{\abs{x}} - x$ is not differentiable at $x = 0$, so smooth first-order optimality conditions are not directly applicable. We now consider two approaches. In the first (\textbf{A1}) we mollify $g(x)$, whereas in the second (\textbf{A2}) we take the variational approach \cite{rockafellar2009variational}.
    \begin{itemize}
        \item[\textbf{A1}.] For $\varepsilon > 0$, consider a Gaussian mollifier $\rho_\varepsilon$. Define the mollified constraint
        \[
        0\geq g_\varepsilon(x) = \int_{\reals} g(x - y)\rho_\varepsilon(y)\d y = \int_{\reals} g(y)\rho_\varepsilon(x - y)\d y.
        \]
        Then, for all $\varepsilon > 0$, we have
        \[
        g_\varepsilon(0) = \int_{\reals} (\sqrt{\abs{y}} - y) \rho_\varepsilon(-y) \d y = \int_{\reals} \sqrt{\abs{y}}\rho_\varepsilon(y) \d y > 0.
        \]
        Thus, $x = 0$ is not a feasible solution for any mollified problem. Moreover, we have
        \[
        g_\varepsilon(x) = \int_{\reals} g(x - y)\rho_\varepsilon(y)\d y
        = \int_{\reals} \sqrt{\abs{x - y}} \rho_\varepsilon(y)\d y - x
        \leq
        \int_{\reals} \sqrt{\abs{y}} \rho_\varepsilon(y)\d y + g(x).
        \]
        Moreover, $\abs{g_\varepsilon(x)-g(x)}\leq\int_{\reals}\sqrt{\abs{y}}\rho_\varepsilon(y)\d y\eqqcolon c_\varepsilon$ for all $x\in\reals$ (since $\abs{\sqrt{\abs{x-y}}-\sqrt{\abs{x}}}\leq\sqrt{\abs{y}}$), and $g_\varepsilon>0$ on $(-\infty,\tfrac{1}{4})$: for $x\leq0$, $g_\varepsilon(x)\geq-x+\int_{\reals}\sqrt{\abs{x-y}}\rho_\varepsilon(y)\d y>0$, whereas for $x\in(0,\tfrac{1}{4})$, averaging $y$ with $-y$ and using $\sqrt{\abs{x-y}}+\sqrt{\abs{x+y}}\geq\sqrt{2\max\{x,\abs{y}\}}\geq\sqrt{x}$ gives $g_\varepsilon(x)\geq\tfrac{1}{2}\sqrt{x}-x>0$. Hence, the feasible set of the mollified problem is contained in $[\tfrac{1}{4},+\infty)\cap\{x\in\reals\st g(x)\leq c_\varepsilon\}$, and the optimizer $x^\ast_\varepsilon$ satisfies $x^\ast_\varepsilon \to 1$ as $\varepsilon \to 0$. That is, the solutions of the mollified problems do not converge to the minimizer of the original problem. In particular, the set of candidate solutions we find with this approach does not include the global minimizer.
        \item[\textbf{A2}.] The regular and general subgradients of $f$ are $\subgradientRegular{f}(0) = \subgradient{f}(0) = \{1\}$. For the constraint set $\feasibleSet = \{x \in \reals \st g(x) \leq 0\}$, we have $\normalCone{\feasibleSet}{0} = \reals$ ($0$ is an isolated feasible point), $\normalCone{\feasibleSet}{1} = (-\infty, 0]$, and $\normalCone{\feasibleSet}{x} = \{0\}$ for $x > 1$. The first-order optimality conditions on a candidate solution $x^\ast$ require that $0 \in \subgradient{f}(x^\ast) + \normalCone{\feasibleSet}{x^\ast}$. Thus, we obtain both the candidate solutions $\{0, 1\}$. Direct inspection, then, reveals the (correct) global optimum $x^\ast = 0$. In particular, the variational approach enables us to find all the local minimizers for the original problem.
    \end{itemize}
}
\end{example}

\subsubsection{Examples of subgradients}
\label{section:ot_wasserstein:functionals}
Next, we study the subgradients of several functionals of interest.
We start with general optimal transport discrepancies, defined in~\eqref{definition:optimaltransport}, which arise, e.g., in uncertainty propagation in dynamical systems \cite{aolaritei2023capture,terpin2023dynamic} and~\gls{acr:dro}~\cite{blanchet2022optimal,kuhn2025distributionally}. With the following result, we generalize~\cite[Proposition 2.17]{Lanzetti2022}.

\includeElement{proposition}{prop_functional_optimal_transport}{true}{proofInMain}

As already observed for the Wasserstein distance, optimal transport discrepancies are subdifferentiable (but not regularly subdifferentiable), and their subgradient can be written in terms of the gradient of the transportation cost and the set of optimal transport plans.
As a notable special case, consider $\varFixed=\delta_{\bar x}$ and $\varReference=\delta_{\hat x}$. Then, $\objective(\varFixed)=\wassersteinDistance{\varFixed}{\varReference}{c}=c(\bar x,\hat x)$ and $\setPlansOptimalCost{\varFixed}{\varReference}{c}=\{\delta_{(\bar x,\hat x)}\}$.
Thus, $\objective$ is differentiable with gradient $\gradient{}\objective(\varFixed)=\pushforward{(\projection{1}{},\gradient{x}c)}{\delta_{(\bar x,\hat x)}}=\delta_{\bar x}\times\delta_{\gradient{x}c(\bar x,\hat x)}$, which is formally analogous to the standard Euclidean gradient of $x\mapsto c(x,\hat x)$.
We now exemplify our results in the case of \WikipediaLink{https://en.wikipedia.org/wiki/Convex_function}{strictly convex} transportation costs:

\includeElement{example}{ex_functional_optimal_transport}{true}{false}

Next, we recall and extend, by characterizing their general subdifferential, the existing subdifferentiability results for expected values in~\cite{ambrosio2005gradient,Lanzetti2022}:

\includeElement{proposition}{prop_functional_expected_value}{true}{proofInMain}

\Cref{proposition:prop_functional_expected_value} characterizes the \emph{general} subdifferential $\subgradient{\objective}$, not only the regular one; as subgradients are transport plans rather than maps, this holds at every $\varFixed\in\Pp{2}{\reals^d}$, including empirical (non-absolutely continuous) measures.
For differentiability, it asks only that $V$ be differentiable, $\lambda$-convex, and bounded above by a quadratic, relaxing~\cite[Proposition 2.20]{Lanzetti2022}, which assumes $V\in\C{2}{\reals^d}{}$ with bounded Hessian---a special case, as any such $V$ is $\lambda$-convex with $\abs{V(x)}\leq A+B\norm{x}^2$.
Thus, the corresponding functional $\objective$ is differentiable and, in particular, $(x,y)\in\supp(\gradient{}\objective)$ if and only if $y=\gradient{} V(x)$. If, instead, $V$ is not bounded below by a quadratic function, then $\objective$ is not proper and might evaluate to $-\infty$ (e.g., $V(x)=-e^{x^2}$ and $\varGeneric=\gaussian{0}{1}$). Similarly, if $V$ grows more than quadratically, then $\objective$ will attain $+\infty$ (e.g., $V(x)=e^{x^2}$ and $\varGeneric=\gaussian{0}{1}$) and cannot be differentiable there.

\includeElement{example}{ex_functional_expected_value}{true}{proofInMain}

Finally, we extend the existing subdifferentiability results for interaction-type functionals~\cite{ambrosio2005gradient,Lanzetti2022}, which model the energy associated with the interaction between particles in a distribution; e.g., see~\cite{santambrogio2015optimal} for an introduction, \cite{buttazzo2012optimal} for an example in physics, and \cite[Section 3.3]{carlier2010optimal} for one in economics.

\includeElement{proposition}{prop_functional_interaction}{true}{proofInMain}

We exemplify the result in the case of quadratic functions:

\includeElement{example}{ex_functional_interaction}{true}{proofInMain}

\subsubsection{Subdifferential calculus}
We now show some useful calculus rules in the Wasserstein space. More can be studied, and we expect this to be a source of interesting future work.

\includeElement{proposition}{prop_calculus}{true}{proofInMain}
This result generalizes \cite[Propositions 2.15 and 2.16]{Lanzetti2022} to our more general subgradients, while relaxing some of the assumptions.
With $g_i(x) = \rho_ix$ for $\rho_1,\rho_2\in\nonnegativeReals$, \cref{proposition:prop_calculus} readily yields the subgradients of the linear combination of functionals:
\includeElement{corollary}{cor_chain_rule_linear}{true}{false}

Next, we deploy~\cref{proposition:prop_calculus} to study the variance via \cref{proposition:prop_functional_expected_value}:
\includeElement{corollary}{cor_functional_variance}{true}{proofInMain}

\subsection{Variational geometry}
\label{section:optimization:var_geom}

\begin{figure}[t]
    \centering
    \begin{minipage}{.32\textwidth}
    \centering
    \includegraphics[width=.8\linewidth]{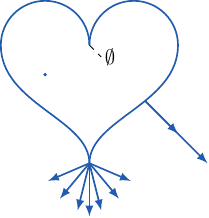}
    \end{minipage}
    \hfill
    \begin{minipage}{.32\textwidth}
    \centering
    \includegraphics[width=.8\linewidth]{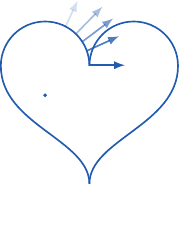}
    \end{minipage}
    \hfill
    \begin{minipage}{.32\textwidth}
    \centering
    \includegraphics[width=.8\linewidth]{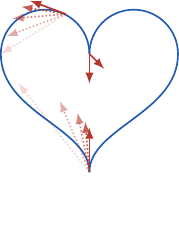}
    \end{minipage}
    \caption{{Variational Geometry in $\reals^d$}. On the left, the arrows depict some examples of regular normals (e.g., at the cusp). The interior points have a trivial normal cone, indicated by a single blue dot. At the same location $x \in \feasibleSet$ they may have different scales and directions. There also may not be any regular normal. In the middle, the \textcolor{ETHblue}{solid blue arrows} depict some examples of normals. They are the limit points of sequences of regular normals. The sequence is colored with increasing opacity as the regular normals approach the general normal. On the right, in \textcolor{ETHred}{dotted red} there are sequences $(x_n-x)/\tau_n \to w \in \tangentCone{\feasibleSet}{x}$, with the tangent vector $w$ being depicted in \textcolor{ETHred}{solid red}. The inner product between tangent vectors and regular normals is nonpositive, whereas the same might not hold with general normals.}
    \label{fig:animation:normal-cones}
\end{figure}

In this section, we study the geometry of the constraint set $\feasibleSet$, which we will combine in~\cref{section:optimization} with the subdifferentiability of the objective function to derive our necessary optimality conditions. Similarly to Euclidean settings, all the necessary information is encoded in tangent and normal cones (see \cref{fig:animation:normal-cones} for an intuition in $\reals^d$). We start with defining the tangent cone:

\includeElement{definition}{def_tangent_cone}{true}{false}

We highlight the formal similarity with the tangent cone in Euclidean spaces, defined by all the $\bar\varTangentReals$ for which there is $(x_n)_{n\in\naturals}\subset\feasibleSet$, $x_n\to\bar x$, and $\tau_n \fromAbove 0$ so that $\frac{1}{\tau_n}(x_n-\bar x)\to \bar\varTangentReals$. In particular, if $\varGeneric_n=\delta_{x_n}$ and $\varFixed=\delta_{\bar x}$, then $\setPlansOptimal{\varFixed}{\varGeneric_n}=\{\delta_{(\bar x,x_n)}\}$ and $\varTangent_n=\delta_{(\bar x, x_n-\bar x)}$. Namely, the definitions in the Wasserstein space and in the Euclidean space are consistent.
Intuitively, the tangent cone is a restriction of the tangent space, and it describes the variations admissible within the set $\feasibleSet$.
Dual to tangent cones are the normal cones:

\includeElement{definition}{def_normal_cone}{true}{false}

Similarly to the definition of subgradient (cf. \cref{definition:def_sub_sup_differentiability}), for general normals we consider convergence in Wasserstein of regular normals. Thus, general normals may not be in the tangent space, and need not be. Analogously to the tangent cone, if we specify the definition of normal cone to Dirac deltas, we recover the notion of normal cone in Euclidean spaces; see also \cref{proposition:prop_normal_cones_trivial}.

Another similarity to Euclidean spaces (cf. \cref{fig:animation:normal-cones}) is the duality between normal and tangent cones:

\includeElement{proposition}{prop_normal_tangent_cones}{true}{proofInMain}

As a first example, we study the normal cones of trivial constraint sets:
\includeElement{proposition}{prop_normal_cones_trivial}{true}{proofInMain}

In particular, analogously to the Euclidean setting, the normal cone at $\varFixed$ trivializes if the feasible set covers the whole space (\cref{proposition:prop_normal_cones_trivial}\ref{item:normal_cone:no_constraints}) or if $\varFixed$ lies in its interior (\cref{proposition:prop_normal_cones_trivial}\ref{item:normal_cone:interior_point}); cf. \cref{fig:animation:normal-cones}. Moreover, when we restrict ourselves to the space of Dirac deltas, the normal cone is consistent with its Euclidean counterpart (\cref{proposition:prop_normal_cones_trivial}\ref{item:normal_cone:delta}); namely, a ``tangent vector'' $\nu$ lies in the normal cone if and only if its ``average'' tangent vector lies in the Euclidean normal cone.
We will study normal cones of level sets of functionals in~\cref{subsec:level_sets}, after having introduced horizon subgradients.

\subsubsection{Technical concepts}
\label{subsec:technical_concepts}
{We now introduce \emph{strong} counterparts of the subgradient,
the normal cone, and the differential, together with a compactness property of
constraint sets. These technical concepts enter (i) the normal-cone characterizations
of \cref{subsec:level_sets}, (ii) the extremal principle and the intersection
rule of \cref{app:intersection}, and (iii) the statement and proof of our
necessary optimality conditions in \cref{section:optimization}. The strong notions are
obtained by requiring the defining inequalities to hold against \emph{all}
transport plans, rather than only the optimal ones:}

\includeElement{definition}{def_strong_sub_sup_differentiability}{true}{false}

{\cref{definition:def_strong_sub_sup_differentiability} differs
from \cref{definition:def_sub_sup_differentiability} only in that the
inequality is required for \emph{all} transport plans
$\setPlans{\varFixed}{\varGeneric}$, and not merely the optimal ones
$\setPlansOptimal{\varFixed}{\varGeneric}$, with the error measured by
$\norm{\varTangent}_{\varFixed}$ (which coincides with
$\wassersteinDistance{\varFixed}{\varGeneric}{2}$ on optimal transport plans).
Analogously, we strengthen the notion of normal cone:}

\includeElement{definition}{def_fuzzy_normal_cone}{true}{false}

{Since membership is now tested against more variations, the
strong notions are more demanding, and the resulting sets are smaller than
their regular counterparts:}

\includeElement{lemma}{strong_vs_regular_normal_cone}{true}{proofInMain}

We strengthen differentiability accordingly:

\includeElement{definition}{def_strong_differentiability}{true}{false}

{By \cref{lemma:strong_vs_regular_normal_cone}\ref{item:strong_vs_regular:subgradient},
applied to $\objective$ and to $-\objective$, a strongly differentiable
functional is differentiable. Hence
\cref{proposition:prop_differentiability_properties}\ref{item:gradient_uniqueness}
applies, and every $\varTangent \in -\subgradientRegularStrong{(-\objective)}(\varFixed)
\cap \subgradientRegularStrong{\objective}(\varFixed)$ satisfies
$\varTangent \localEquiv{\varFixed} \gradient{}{\objective}(\varFixed)$; we call
this element the \emph{strong gradient} $\gradientStrong{}{\objective}(\varFixed)$
(so $\gradientStrong{}{\objective}(\varFixed) \localEquiv{\varFixed} \gradient{}{\objective}(\varFixed)$).}

{A few comments are in order. First, the distinction between
the strong and the standard notions is specific to the Wasserstein space:
There is exactly one way to couple a Dirac delta $\delta_{\bar x}$ with any
probability measure (namely, their product), so at Dirac deltas the strong
regular subgradient and the strong regular normal cone coincide with their
counterparts in
\cref{definition:def_sub_sup_differentiability,definition:def_normal_cone}.
Accordingly, no analogous distinction arises in Euclidean spaces. Second, the
interplay between the two families is dictated by the direction of our
results. Our optimality conditions are \emph{necessary}: They confine the
candidate minimizers to sets built from subgradients and normal cones (cf.
\cref{theorem:th_first_order_optimality_conditions}), and they remain valid
if any of these sets is replaced by a superset---the conditions merely ``get
bigger'', yielding more candidates but never missing a minimizer. Thus,
whenever an exact characterization of a normal cone is out of reach, an
upper bound suffices for the purpose of necessary conditions, as in
\cref{subsec:level_sets}. Conversely, since smaller sets make a necessary
condition sharper, we state the conclusion of
\cref{theorem:th_first_order_optimality_conditions} with the \emph{strong}
normal cone; the version with the normal cone follows from
\cref{lemma:strong_vs_regular_normal_cone}\ref{item:strong_vs_regular:normal}.
Assumptions, by contrast, become \emph{weaker} when imposed on the
smaller, strong sets: A requirement on all strong regular normals is easier
to satisfy than the same requirement on all regular normals, and the
corresponding results apply more broadly. We exploit this asymmetry three
times: in the constraint qualifications of the intersection rule
(\cref{app:intersection}) and of our first-order optimality conditions
(\cref{theorem:th_first_order_optimality_conditions}), which feature the
strong normal cones, and in the compactness property below, which is imposed
only along strong regular normals.}

{Finally, to prove our first-order optimality conditions in~\cref{theorem:th_first_order_optimality_conditions}, we will require a technical property for the constraint set:
\includeElement{definition}{def_snc}{true}{false}}

{Informally, \cref{definition:def_snc} asks that, along
sequences in $\feasibleSet$, narrow convergence of strong regular normals
upgrades to convergence in Wasserstein. In finite dimensions, where weak and
strong convergence coincide, every set enjoys the analogous property, and no
such condition appears in Euclidean variational
analysis~\cite{rockafellar2009variational,mordukhovich2018variational}; in
infinite dimensions, it is a standard regularity
requirement~\cite{mordukhovich2006variationalI}. The same is true in the
Wasserstein space, where narrow convergence does not imply convergence in
Wasserstein.
In line with the discussion above, imposing
this property only along the \emph{strong} regular normals yields a weaker
assumption and, in turn, more widely applicable optimality conditions. We
verify the \gls*{acr:snc} property for closed Wasserstein balls in
\cref{proposition:prop_normal_cone_wasserstein_ball}.}

\subsubsection{Variational geometry and epigraphs}
\label{section:variational-geometry:epigraphs}
In \cref{section:optimization}, we provide the first-order optimality conditions for general, constrained, optimization problems and present differentiable functionals as a special case.
From a technical standpoint, however, we \emph{first} establish optimality conditions for differentiable functions, and \emph{then} study the nondifferentiable case via an epigraphical argument. This will yield a ``constraint qualification'', which is typical in Euclidean spaces and we find in the Wasserstein space as well.
More formally, recall that the \WikipediaLink{https://en.wikipedia.org/wiki/Epigraph_(mathematics)}{epigraph} of a proper functional $\objective: \Pp{2}{\reals^d} \to \realsBar$ is the set
$\widetilde{\mathrm{epi}}\objective
\coloneqq
\left\{
(\varGeneric, \beta)
\st
\varGeneric \in \Pp{2}{\reals^d}, \objective(\varGeneric) \leq \beta
\right\}
\subset \Pp{2}{\reals^d}\times\reals$.
Then, we immediately observe that
$\inf_{\varGeneric \in \feasibleSet}\objective(\varGeneric)
=
\inf_{(\varGeneric, \beta) \in \widetilde{\mathrm{epi}}(\objective) \cap (\feasibleSet\times\reals)} \beta
=
\inf_{(\varGeneric, \beta) \in \widetilde{\mathrm{epi}}(\objective) \cap (\feasibleSet\times\reals)} \tilde{l}(\varGeneric, \beta)$,
with $\tilde{l}: \Pp{2}{\reals^d}\times\reals \to \reals$ being the linear function $\tilde{l}(\varGeneric, \beta) = \beta$. Being linear, $\tilde{l}$ is in principle easy to optimize.

Unfortunately, this argument does not directly carry over to the Wasserstein space: While on $\reals^d$ the epigraph is a subset of $\reals^{d+1}$, on $\Pp{2}{\reals^d}$, it is a subset of $\Pp{2}{\reals^d}\times\reals$ and not $\Pp{2}{\reals^{d+1}}$.
Thus, one has to study the differential structure and variational geometry of $\Pp{2}{\reals^d}\times\reals$.
To avoid this burden, we take a different angle. Since each $x\in\reals^d$ can be ``embedded'' in the space of probability measures as the Dirac delta $\delta_x\in\Pp{2}{\reals^d}$,
we define the epigraph as
{\begin{equation*}
\epigraph{\objective}
\coloneqq
\left\{
\varGeneric\times\delta_\beta
\st
\varGeneric \in \Pp{2}{\reals^d}, \objective(\varGeneric) \leq \beta
\right\} \subseteq \Pp{2}{\reals^{d+1}}.
\end{equation*}}
In particular, contrary to $\widetilde{\mathrm{epi}}\objective$, $\epigraph{\objective}$ is now a subset of $\Pp{2}{\reals^{d+1}}$.
Then, the optimization of $\objective$ over $\feasibleSet$ is equivalent to the optimization of $l:\Pp{2}{\reals^{d+1}} \to \reals$, $l(\boldsymbol{\varGeneric}) = \expectedValue{\boldsymbol{x} \sim \boldsymbol{\varGeneric}}{\innerProduct{(0, \ldots, 0, 1)}{\boldsymbol{x}}}$, over $\epigraph{\objective} \cap \feasibleSet_e$, with
$
\feasibleSet_e \coloneqq \{
\varGeneric \times \delta_\beta
\st
\varGeneric \in \feasibleSet, \beta \in \reals
\}.
$
{Clearly, the same equivalence holds locally, which is what we will use to
derive optimality conditions.}
Note that we use bold symbols when working in $\reals^{d+1}$ and $\Pp{2}{\reals^{d+1}}$ (e.g., $\boldsymbol{x}=(x,\beta)$ with $x\in\reals^d$ and $\beta\in\reals$).
Given $\boldsymbol{\varGeneric}$, we can ``extract'' $\varGeneric$ and $\delta_\beta$ via $\pushforward{(\projection{1}{},\ldots,\projection{d}{})}{\boldsymbol{\varGeneric}}=\varGeneric$ (i.e., we forget the last component) and $\pushforward{(\projection{d+1}{})}{\boldsymbol{\varGeneric}}=\delta_\beta$ (i.e., we keep only the last component).
Similarly, given a variation $\boldsymbol{\varTangent}\in\Pp{2}{\reals^{d+1}\times\reals^{d+1}}$ such that $\pushforward{\projection{1}{}}{\boldsymbol{\varTangent}} = \boldsymbol{\varGeneric}$, we can extract the variation corresponding to $\varGeneric$ via $\pushforward{T}{\boldsymbol{\varTangent}}$, where
\begin{equation}
\label{equation:extraction-map-epigraph}
T\coloneqq (\projection{1}{},\ldots,\projection{d}{},\projection{d+2}{},\ldots,\projection{2d+1}{})
\end{equation}
(i.e., we forget the last component of the first and the second marginal),
and the variation ``corresponding'' to $\delta_\beta$ via $\pushforward{(\projection{d+1}{},\projection{2d+2}{})}{\boldsymbol{\varTangent}}$ (i.e., we only keep the last component of the first and the second marginal).
With this formalism, the gradient of $l$ follows directly from \cref{proposition:prop_functional_expected_value} and it reads $\boldsymbol{\varGeneric}\times\boldsymbol{\delta_{(0,\ldots,0,1)}}$. Thus, the complexity of our main result moves to proving first-order optimality conditions in the differentiable case subject to epigraphical constraints, for which, in particular, we need to study the variational geometry of $\epigraph{\objective}$, which we do next.
To start, observe that epigraphical normals ``point downwards in expectation'':

\includeElement{proposition}{prop_epigraph_normal_tangent}{true}{proofInMain}

\cref{proposition:prop_epigraph_normal_tangent} is reminiscent of the well-known fact that, in Euclidean settings, epigraphical normals never point upwards.
We depict this idea below on the right.

\vspace{0.5\baselineskip}

\begin{center}
\begin{minipage}{.65\linewidth}
For the (discontinuous) function $f$ with graph depicted in black, the epigraph is the region in light green above the graph. The thick blue arrows represent the pairs $(\varTangentReals, -1)$, with $\varTangentReals$ regular subgradient of $f$ at the discontinuity point (black-filled point). These arrows are also related to the normals of the epigraph, which contain their $\lambda$-scaled versions as well. These are depicted in dashed blue, and they intersect the $-\lambda$ plane.
\end{minipage}
\hfill
\begin{minipage}{.32\linewidth}
\includegraphics[width=.9\linewidth]{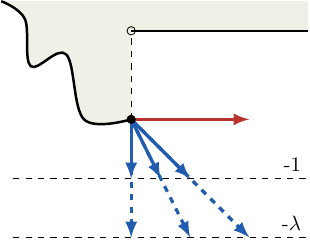}
\end{minipage}
\end{center}

\vspace{\baselineskip}

We rigorously show that the relationship between the nonhorizontal normals (in expectation) to $\epigraph{\objective}$ and the subgradients of $\objective$ translates to the Wasserstein space:
\includeElement{proposition}{prop_subgradients_epigraphical_normals}{true}{proofInMain}
Note that \cref{proposition:prop_subgradients_epigraphical_normals} characterizes \emph{all} nonhorizontal normals.
Indeed, if $\boldsymbol{\varTangent}\in \normalCone{\epigraph{\objective}}{\varFixed\times\delta_{\objective(\varFixed)}}$ satisfies $\expectedValue{(\boldsymbol{x},\boldsymbol{\varTangentReals})\sim\boldsymbol{\varTangent}}{\innerProduct{(0,\ldots,0,1)}{\boldsymbol{\varTangentReals}}}=-\lambda<0$, then, since the normal cone is a cone (cf. \cref{proposition:prop_normal_tangent_cones}), $\frac{1}{\lambda}\boldsymbol{\varTangent}$ also belongs to $\normalCone{\epigraph{\objective}}{\varFixed\times\delta_{\objective(\varFixed)}}$ and it satisfies $\expectedValue{(\boldsymbol{x},\boldsymbol{\varTangentReals})\sim\frac{1}{\lambda}\boldsymbol{\varTangent}}{\innerProduct{(0,\ldots,0,1)}{\boldsymbol{\varTangentReals}}}=
\expectedValue{(\boldsymbol{x},\boldsymbol{\varTangentReals})\sim\boldsymbol{\varTangent}}{\innerProduct{(0,\ldots,0,1)}{\frac{1}{\lambda}\boldsymbol{\varTangentReals}}}=
-1$. Thus,
$\pushforward{T}{\frac{1}{\lambda}\boldsymbol{\varTangent}}\in\subgradient{\objective}(\varFixed)$.
In Euclidean spaces, horizontal normals are referred to as horizon subgradients (we depict one in red in the figure above).
Thus, we introduce the natural counterpart in the Wasserstein space:
\includeElement{definition}{def_sub_sup_differentiability_horizon}{true}{proofInMain}

The horizon subgradients play a role in the constraint qualifications and, as we shall see in~\cref{section:optimization}, are required for a comprehensive theory of optimization. Similar ideas are necessary already in Euclidean spaces and, perhaps surprisingly, they naturally carry over to the Wasserstein space. Since they are defined starting from the general normals to the epigraph, they are not guaranteed to be in the tangent space.
Nonetheless, for sufficiently well-behaved functionals they are trivial:

\includeElement{proposition}{prop_subgradients_horizon_trivial}{true}{proofInMain}

For instance, all locally Lipschitz functionals are strictly continuous and therefore have trivial horizon subgradients.
This is also the case for many of the functionals discussed so far:
\includeElement{proposition}{cor_subgradients_horizon}{true}{proofInMain}

The strict-continuity assumption in \cref{proposition:prop_subgradients_horizon_trivial} cannot be weakened to differentiability: a functional may be differentiable at $\varFixed$ and still have a nontrivial horizon subgradient there.
\begin{example}[A differentiable functional with a nontrivial horizon subgradient]
\label{example:differentiable_nontrivial_horizon}
Let $d=1$ and let $f\colon\reals\to\reals$ be $f(0)\coloneqq 0$ and $f(x)\coloneqq x^2\sin(x^{-2})$ for $x\neq 0$. Then $f$ is differentiable on all of $\reals$, with $f'(0)=0$ and, for $x\neq 0$, $f'(x)=2x\sin(x^{-2})-2x^{-1}\cos(x^{-2})$; in particular, $f'$ is unbounded near $0$ and takes arbitrarily large values of both signs. Consider $\objective\colon\Pp{2}{\reals}\to\reals$, $\objective(\varGeneric)\coloneqq f\bigl(\expectedValue{x\sim\varGeneric}{x}\bigr)$, and $\varFixed\coloneqq\delta_0$, so that $\objective(\delta_0)=0$.

\emph{$\objective$ is differentiable at $\delta_0$, with $\gradient{}{\objective}(\delta_0)\localEquiv{\delta_0}\localZero{\delta_0}$.} For $\varGeneric\in\Pp{2}{\reals}$ and $\varTangent\in\pushforward{(\projection{1}{},\projection{2}{}-\projection{1}{})}{\setPlansOptimal{\delta_0}{\varGeneric}}$, the mean of $\varGeneric$ equals $\expectedValue{(x,\varTangentReals)\sim\varTangent}{\varTangentReals}$, so $\abs{\expectedValue{x\sim\varGeneric}{x}}\leq\norm{\varTangent}_{\delta_0}=\wassersteinDistance{\delta_0}{\varGeneric}{2}$. Since $f$ is differentiable at $0$ with $f'(0)=0$, we get $\objective(\varGeneric)-\objective(\delta_0)=f(\expectedValue{x\sim\varGeneric}{x})-f(0)=\onotation{\abs{\expectedValue{x\sim\varGeneric}{x}}}=\onotation{\norm{\varTangent}_{\delta_0}}$. Hence $\localZero{\delta_0}$ is a regular subgradient of both $\objective$ and $-\objective$ at $\delta_0$, and \cref{definition:def_differentiability} gives differentiability with $\gradient{}{\objective}(\delta_0)\localEquiv{\delta_0}\localZero{\delta_0}$.

\emph{$\objective$ is not strictly continuous at $\delta_0$.} Indeed, $\objective(\delta_s)=f(s)$ and $\wassersteinDistance{\delta_{s_1}}{\delta_{s_2}}{2}=\abs{s_1-s_2}$, so $\frac{\abs{\objective(\delta_{s_1})-\objective(\delta_{s_2})}}{\wassersteinDistance{\delta_{s_1}}{\delta_{s_2}}{2}}=\frac{\abs{f(s_1)-f(s_2)}}{\abs{s_1-s_2}}$, whose supremum as $s_1,s_2\to 0$ is infinite because $f'$ is unbounded near $0$.

\emph{The horizon subgradient of $\objective$ at $\delta_0$ is nontrivial.} Fix $w\neq 0$ and choose $x_n\to 0$ with $f'(x_n)\to+\infty$ if $w>0$ (resp., $f'(x_n)\to-\infty$ if $w<0$), and set $\tau_n\coloneqq w/f'(x_n)\fromAbove 0$. Repeating the estimate above at $x_n$ (against epigraph points $\varGeneric\times\delta_\beta$ with $\beta\geq\objective(\varGeneric)$) shows that the deterministic variation $\boldsymbol{\varTangent}_n\coloneqq\delta_{\left((x_n,f(x_n)),\,(f'(x_n),-1)\right)}$ belongs to $\normalConeRegular{\epigraph{\objective}}{\delta_{x_n}\times\delta_{f(x_n)}}$, with last velocity component $-1$ in expectation and $\pushforward{T}{\boldsymbol{\varTangent}_n}=\delta_{(x_n,f'(x_n))}$ ($T$ as in \eqref{equation:extraction-map-epigraph}). Since the regular normal cone is a cone, $\tau_n\boldsymbol{\varTangent}_n=\delta_{\left((x_n,f(x_n)),\,(w,-\tau_n)\right)}$ is also a regular normal. As $n\to\infty$, the base points $\delta_{x_n}\times\delta_{f(x_n)}\in\epigraph{\objective}$ converge to $\delta_0\times\delta_0$ and $\tau_n\boldsymbol{\varTangent}_n\convergenceWasserstein\boldsymbol{\varTangent}\coloneqq\delta_{\left((0,0),(w,0)\right)}$, so $\boldsymbol{\varTangent}\in\normalCone{\epigraph{\objective}}{\delta_0\times\delta_0}$ by \cref{definition:def_normal_cone}. Its last velocity component has zero expectation and $\pushforward{T}{\boldsymbol{\varTangent}}=\delta_{(0,w)}=\pushforward{(\identity{},w)}{\delta_0}$, so $\pushforward{(\identity{},w)}{\delta_0}\in\subgradientHzn{\objective}(\delta_0)$ by \cref{definition:def_sub_sup_differentiability_horizon}. As $w\neq 0$ was arbitrary, $\subgradientHzn{\objective}(\delta_0)\neq\{\localZero{\delta_0}\}$.

{Thus $\objective$ is differentiable at $\delta_0$ yet has a nontrivial horizon subgradient, so the hypothesis of \cref{proposition:prop_subgradients_horizon_trivial} is sharp. Consistently, when we treat the differentiable case of \cref{theorem:th_first_order_optimality_conditions} we do not rely on the horizon subgradient being trivial.
}
\end{example}

With horizon subgradients, we can finally provide a full characterization of epigraphical normals:

\includeElement{proposition}{th_subgradients_epigraphical_normals}{true}{proofInMain}

\paragraph{Discussion.}
The definition of horizon subgradients introduced in \cref{definition:def_sub_sup_differentiability_horizon} follows the geometric ideas used for Euclidean spaces in \cite[Definition 1.18]{mordukhovich2018variational}. It departs from the alternative definition of the horizon subgradients as the scaled limit of regular subgradients \cite[Definition 8.3]{rockafellar2009variational}, according to which $\varTangentFixed \in \subgradientHzn{\objective}(\varFixed)$ if there are sequences $\varGeneric_n \to \varFixed, \objective(\varGeneric_n) \to \objective(\varFixed)$ and $\varTangent_n \in \subgradientRegular{\objective}(\varGeneric_n)$ so that there exists $(\tau_n)_{n \in \naturals} \subseteq \nonnegativeReals, \tau_n \fromAbove 0$ and $\tau_n\varTangent_n \convergenceWasserstein \varTangentFixed$.
Nonetheless, it is well known that, at least in Euclidean settings, the two definitions are identical, cf. \cite[Theorem 8.9]{rockafellar2009variational}.
Finally, we could have defined the subgradients also as the nonhorizontal epigraphical normals, as in~\cite[Definition 1.18]{mordukhovich2018variational}. Nonetheless, to ease the presentation, we opted for the more established~\cref{definition:def_sub_sup_differentiability}, as in~\cite[Definition 8.3]{rockafellar2009variational}.

\subsubsection{Example: Level sets}\label{subsec:level_sets}

To conclude this section, we study the normal cone to sublevel sets of continuous functionals:
\includeElement{proposition}{prop_normal_cones_constraints_simple}{true}{proofInMain}
\cref{proposition:prop_normal_cones_constraints_simple} is closely related to Lagrange multipliers in the Wasserstein space: The normal cone to a sublevel set of a functional $\objective{}$ consists of any scaling (by a multiplier) of the subgradients of $\objective{}$. In particular, it includes the case where the constraint is not active (i.e., the multiplier is zero and the normal cone trivializes) and the case where the constraint is active (i.e., the multiplier is nonzero). The condition $\localZero{\varFixed} \not\in \subgradient{\objective}(\varFixed)$ resembles the same ``constraint qualification'' required in Euclidean settings (cf. \cite[Proposition 10.3]{rockafellar2009variational}), and in practice it is often not restrictive.
With this result, we now study the normal cone to (closed) Wasserstein balls:
\includeElement{proposition}{prop_normal_cone_wasserstein_ball}{true}{proofInMain}
As a special case, we study the feasible set of the optimization problem \eqref{equation:problem:intro:simple}, presented in \cref{sec:introduction}, namely the second-moment constraint:
\begin{example}[Second-moment constraint]
\label{example:ex_second_moment_constraint}
For $\varepsilon > 0$, consider
\begin{equation*}
\feasibleSet = \{\varGeneric \in \Pp{2}{\reals^d}\st \expectedValue{x \sim \varGeneric}{\norm{x}^2} \leq \varepsilon^2\}.
\end{equation*}
Then, $\feasibleSet = \closedWassersteinBall{\delta_0}{\varepsilon}{2}$ and, thus, for any $\varFixed \in \Pp{2}{\reals^d}$ in the interior of $\feasibleSet$ the (strong) normal cone trivializes, whereas for all the ones with
$
\expectedValue{x \sim \varFixed}{\norm{x}^2}= \varepsilon^2
$
we have
$
\normalConeStrong{\feasibleSet}{\varFixed}
\subseteq
\left\{\pushforward{(\identity{}, 2\lambda\identity{})}{\varFixed} \st \lambda \geq 0\right\}.
$
\end{example}

\section{Optimality conditions}\label{section:optimization}
In this section, we present and discuss our main result. To start, we need a notion of local optimality:
\includeElement{definition}{def_wasserstein_space_local_optimum}{true}{false}

Next is the main result of this work, which encompasses the celebrated Karush--Kuhn--Tucker and Lagrange conditions: first-order necessary optimality conditions in the Wasserstein space.
\includeElement{theorem}{th_first_order_optimality_conditions}{true}{proofInMain}
Geometrically, the constraint qualification requires that the feasible set $\feasibleSet$ is not ``tangent'' to the epigraph at $\varOptimal$.
Similarly to the Euclidean setting, at optimality at least one element of the subgradient of the cost functional (direction of improvement) must align and have the same magnitude (but opposite direction) with one element of the normal cone (the ``blocked'' directions). We refer to \cref{section:optimization:discussion} for further discussion.
In the unconstrained case, the normal cone trivializes (cf.~\cref{proposition:prop_normal_cones_trivial}) and we recover Fermat's rule:
\includeElement{theorem}{th_fermat_rule}{true}{proofInMain}

\subsection{Discussion}
\label{section:optimization:discussion}

Before showcasing our first-order necessary optimality conditions in various examples, we discuss several related aspects.

\paragraph{Relation to classical variational analysis.}
The conditions in~\cref{theorem:th_first_order_optimality_conditions,theorem:th_fermat_rule} very much resemble their counterparts in Euclidean spaces and Asplund spaces: the negative subgradient of the objective function must lie in the normal cone to the feasible set, provided that a constraint qualification (i.e., the intersection of the horizon subgradient and the ``negative'' normal cone is trivial) and a regularity condition on the feasible set hold (i.e., the set $\feasibleSet$ is \gls*{acr:snc}); e.g., see~\cite{rockafellar2009variational,mordukhovich2018variational} for the corresponding statement in Euclidean spaces (where the \gls*{acr:snc} condition always holds and thus does not appear) and \cite[Section 5.1]{mordukhovich2006variationalII} for the infinite-dimensional version.
Note that these results in Euclidean spaces generalize, among others, the celebrated Karush--Kuhn--Tucker and Lagrange conditions.

\paragraph{Comparison with the literature.}
General necessary optimality conditions in the Wasserstein space are studied in~\cite{Lanzetti2022}, in the smooth (or nonsmooth but convex) setting where the constraint is the level set of a real-valued functional, and in~\cite{bonnet2019pontryagin2,bonnet2021necessary,bonnet2019pontryagin} in the different setting of optimal control problems in the Wasserstein space. With this work, we provide a general theory of variational analysis in the Wasserstein space and significantly generalize existing necessary conditions for static optimization.

\paragraph{How to apply~\cref{theorem:th_first_order_optimality_conditions}.}
{
The application of~\cref{theorem:th_first_order_optimality_conditions} requires characterizing the subgradients of $\objective$ and the normal cone to $\feasibleSet$.
For many functionals of practical interest, these characterizations are provided in~\cref{section:ot_wasserstein:functionals}. Moreover,~\cref{proposition:prop_calculus} provides tools to derive the subgradients of additional functionals.
The characterization of subgradients for other functionals is left to future research. Importantly, this needs to be carried out only once for each functional and, as in Euclidean spaces, we envision a growing collection of reusable results.
Similarly,~\cref{proposition:prop_normal_cones_constraints_simple} provides a formula for the normal cones to the important class of feasible sets defined as level sets of functionals. These formulas can then be combined using the product and intersection rule in~\cref{app:intersection}. As with subgradients, the derivation of further normal-cone formulas is left to future research.
Finally, the only remaining assumptions are that (i) the constraint set is \gls*{acr:snc} (properness is typically immediate in meaningful optimization settings) and (ii) the constraint qualification $\subgradientHzn{\objective}(\varOptimal) \cap\left(-\normalConeStrong{\feasibleSet}{\varOptimal} \right) = \{\localZero{\varOptimal}\}$, imposed only along the smaller strong normal cone (cf. \cref{subsec:technical_concepts}).
The \gls*{acr:snc} assumption is not required when the objective function is differentiable or when the problem is unconstrained. If the objective is merely subdifferentiable, then one must instead study the regularity properties of the constraint set, as we did in~\cref{proposition:prop_normal_cone_wasserstein_ball} for Wasserstein balls.
By contrast, the constraint qualification is automatically satisfied for many functionals of practical interest, including nonsmooth ones, because their horizon subdifferentials are trivial. When the horizon subdifferential is nontrivial, however, this condition must be verified explicitly. In the unconstrained case, the normal cone reduces to $\{\localZero{\varOptimal}\}$ and the constraint qualification holds automatically.
}

\paragraph{Computational aspects.}
Our optimality conditions transform the problem of minimizing a function over the infinite-dimensional space of probability measures into an inclusion problem at the level of transport plans which, in turn, amounts to a set of conditions over $\reals^d\times\reals^d$. Depending on the application, they can be solved in (quasi) closed-form or addressed via (nonlinear) function approximators \cite{amos2017input,uscidda2023monge}, as we show in our examples. For this reason, we argue that our necessary conditions are significantly more tractable than the initial infinite-dimensional optimization problem.
Nonetheless, more rigorous statements of the computational aspects of necessary conditions for optimality depend on the specifics of the problem at hand, as is already the case in Euclidean settings.

\paragraph{Sufficient conditions.}
As in Euclidean settings, our conditions are not sufficient for optimality.
We expect sufficient conditions for optimality to be intimately related to geodesic convexity~\cite[Section 7]{ambrosio2005gradient} and second-order calculus in the Wasserstein space~\cite{gigli2012second}; see~\cite{Lanzetti2022} for preliminary results.

\subsection{Pedagogical examples}
In the remainder of this section, we illustrate our necessary conditions across several examples. Despite its pedagogical approach, this section culminates with \cref{example:mean-var}, which already generalizes existing insights in \gls*{acr:dro}.
We start by rigorously solving the example in \cref{sec:introduction}:
\begin{example}[Expected value subject to second-moment constraints]
\label{example:max:exp:second-moment}
For $\theta \neq 0$ and $\varepsilon > 0$, consider the problem
    $\inf_{\varGeneric \in \Pp{2}{\reals^d}}\; \objective(\varGeneric)$, where $\objective(\varGeneric)\coloneqq\expectedValue{x \sim \varGeneric}{\innerProduct{\theta}{x}}$, subject to
    $\expectedValue{x \sim \varGeneric}{\norm{x}^2} \leq \varepsilon^2$.
By \cref{example:ex_functional_expected_value}, $\objective$ is differentiable with gradient $\gradient{}{\objective}(\varGeneric) = \pushforward{(\identity{}, \theta)}{\varGeneric}$ and general subdifferential $\subgradient{\objective}(\varGeneric)=\{\gradient{}{\objective}(\varGeneric)\}$; its horizon subgradient is trivial (\cref{proposition:cor_subgradients_horizon}), so the constraint qualification of \cref{theorem:th_first_order_optimality_conditions} holds.
The strong normal cone to the second-moment constraint was studied in \cref{example:ex_second_moment_constraint}, $\normalConeStrong{\feasibleSet}{\varFixed} \localEquiv{\varFixed} \{\localZero{\varFixed}\}$ if $\expectedValue{x \sim \varFixed}{\norm{x}^2} < \varepsilon^2$ and $\normalConeStrong{\feasibleSet}{\varFixed} \subseteq \{\pushforward{(\identity{}, 2\lambda\identity{})}{\varFixed} \st \lambda \geq 0\}$ if $\expectedValue{x \sim \varFixed}{\norm{x}^2} = \varepsilon^2$; the set $\feasibleSet$ is a closed Wasserstein ball, hence closed, and it is \gls*{acr:snc} by \cref{proposition:prop_normal_cone_wasserstein_ball}. Any candidate optimal solution $\varOptimal \in \Pp{2}{\reals^d}$ must satisfy the optimality conditions in \cref{theorem:th_first_order_optimality_conditions}.
We study two cases. If $\expectedValue{x \sim \varOptimal}{\norm{x}^2} < \varepsilon^2$, it holds
$\localZero{\varOptimal}
\localEquiv{\varOptimal}
\gradient{}{\objective}(\varOptimal)
$ or, equivalently,
$0 = \norm{\gradient{}{\objective}(\varOptimal)}_{\varOptimal}^2
=\int_{\reals^d} \norm{\theta}^2 \d\varOptimal(x) > 0$,
which is absurd. Therefore, the optimal solution (if it exists) lies at the boundary, and a $\lambda > 0$ (for $\lambda = 0$, the condition reduces to the interior case, which is absurd) and a ``sum coupling'' $\varTangentCoupling \in \setPlansCommonMarginal{\gradient{}{\objective}(\varOptimal)}{\lambda \pushforward{(\identity{}, 2\identity{})}{\varOptimal}}{1}$ must exist so that
$\localZero{\varOptimal}
\localEquiv{\varOptimal}
\gradient{}{\objective}(\varOptimal) \localSum{\varTangentCoupling}
\lambda \pushforward{(\identity{}, 2\identity{})}{\varOptimal}$ or, equivalently,
$0 = \norm{\gradient{}{\objective}(\varOptimal) \localSum{\varTangentCoupling}
\lambda\pushforward{(\identity{}, 2\identity{})}{\varOptimal}}_{\varOptimal}^2
=\int_{\reals^d} \norm{\theta + 2\lambda x}^2 \d\varOptimal(x)$.
The above is precisely a ``Lagrange's condition'' in the Wasserstein space, which yields $\varOptimal = \delta_{-\frac{\theta}{2\lambda}}$ and $\lambda = \frac{\norm{\theta}}{2\varepsilon}$, since $\expectedValue{x \sim \varOptimal}{\norm{x}^2} = \varepsilon^2$.
\end{example}

We now show that necessary conditions for optimality can be used to produce certificates that an optimal solution does not exist:

\begin{example}[A certificate of nonexistence of a minimizer]
\label{example:ex_infimum_not_attained}
Consider the problem of finding the worst-case mean-variance risk over the full Wasserstein space. In this case, for $\rho > 0$ and $0\neq \theta\in\reals^d$, we seek the minimum of the functional
\begin{equation*}
    \objective(\varGeneric) \coloneqq
    -\left(
    \expectedValue{x \sim \varGeneric}{\innerProduct{\theta}{x}}
    +
    \frac{\rho}{2}\variance{x \sim \varGeneric}{\innerProduct{\theta}{x}}
    \right).
\end{equation*}
By~\cref{proposition:prop_functional_expected_value,proposition:prop_calculus,corollary:cor_functional_variance}, $\objective$ is differentiable with gradient
    $\gradient{}{\objective}(\varGeneric) \localEquiv{\varFixed} -\pushforward{(\identity{}, (1 + \rho(\innerProduct{\theta}{\identity{}} - \expectedValue{x \sim \varGeneric}{\innerProduct{\theta}{x}}))\theta)}{\varGeneric}$.
To find the unconstrained minimizer of $\objective$, we can deploy \cref{theorem:th_fermat_rule}:
    {at any optimal $\varOptimal \in \Pp{2}{\reals^d}$ it holds $\localZero{\varOptimal} \localEquiv{\varOptimal}\gradient{}{\objective}(\varOptimal)$ or, equivalently, $0=\norm{\gradient{}{\objective}(\varOptimal)}^2_{\varOptimal}
    =\int_{\reals^d}\norm{(1 + \rho(\innerProduct{\theta}{x} - \expectedValue{y \sim \varOptimal}{\innerProduct{\theta}{y}}))\theta}^2\d \varOptimal(x)$.
Thus, }
\begin{equation}
\label{eq:contr}
    0 =1 + \rho(\innerProduct{\theta}{x} - \expectedValue{y \sim \varOptimal}{\innerProduct{\theta}{y}}) \quad\almostEverywhere{\varOptimal},
\end{equation}
but taking the expected value with respect to $\varOptimal$ of both sides in \eqref{eq:contr} yields (it holds $\almostEverywhere{\varOptimal}$)
$0 = \expectedValue{x \sim \varOptimal}{0} = \expectedValue{x \sim \varOptimal}{1 + \rho(\innerProduct{\theta}{x} - \expectedValue{y \sim \varOptimal}{\innerProduct{\theta}{y}})} = 1$,
a contradiction. Thus, the infimum of $\objective$ is not attained.
\end{example}
While not a novel result, the above example shows that with only simple algebraic manipulations we can use the optimality conditions in \cref{theorem:th_first_order_optimality_conditions,theorem:th_fermat_rule} to prove nonexistence of a solution, in a way that mirrors arguments in Euclidean spaces.

We now consider an example in the setting of~\gls*{acr:dro}, a ubiquitous framework for decision-making under uncertainty. We seek to evaluate the worst-case mean-variance of a linear portfolio with allocation $\theta$ but now constrained to a closed Wasserstein ball of a given radius $\varepsilon$ and center $\varReference$.
Formally, this example is similar to~\cref{example:ex_infimum_not_attained}, but now the worst-case is taken over a Wasserstein ball instead of the entire Wasserstein space.
As we shall see below, our necessary conditions lead to a closed-form solution for the worst-case probability measure, generalizing~\cite[Section 4]{Lanzetti2022} to non-absolutely continuous measures (which, among others, allows us to study the data-driven setting where $\varReference$ is empirical) and~\cite[Appendix E]{nguyen2021mean}, which does not provide a closed-form solution and assumes positive variance of $\varReference$.

\begin{example}[A constrained problem: mean-variance \gls*{acr:dro}]
\label{example:mean-var}
Consider the mean-variance functional $\objective: \Pp{2}{\reals^d}\to\realsBar$ defined in \cref{example:ex_infimum_not_attained} and the constraint $\varGeneric \in \closedWassersteinBall{\varReference}{\varepsilon}{2}$ for some $\varReference \in \Pp{2}{\reals^d}$ and $\varepsilon>0$.
We computed the gradient of $\objective$ in~\cref{example:ex_infimum_not_attained}, and we recall from~\cref{proposition:prop_normal_cone_wasserstein_ball} that the strong normal cone to $\closedWassersteinBall{\varReference}{\varepsilon}{2}$ at any $\varFixed \in \closedWassersteinBall{\varReference}{\varepsilon}{2}$ satisfies
$\normalConeStrong{\closedWassersteinBall{\varReference}{\varepsilon}{2}}{\varFixed}
\subseteq
\{
\varTangent\st
\varTangent\in
\pushforward{(\projection{1}{}, 2\lambda(\projection{1}{} - \projection{2}{}))}{\setPlansOptimal{\varFixed}{\varReference}},
\lambda\geq 0
\}$, and that the (closed) ball is \gls*{acr:snc}.
{The remaining hypotheses of
\cref{theorem:th_first_order_optimality_conditions} hold as well: $\objective$
is a linear combination of an expected value and a variance of the linear
function $x \mapsto \innerProduct{\theta}{x}$ (whose square has uniformly
bounded Hessian $2\theta\transpose{\theta}$), so each term is Lipschitz on
bounded sets (\cref{proposition:cor_subgradients_horizon}), $\objective$ is
strictly continuous, and its horizon subgradient is trivial
(\cref{proposition:prop_subgradients_horizon_trivial}), whence the constraint
qualification; moreover, $\subgradient{\objective}(\varOptimal) =
\{\gradient{}{\objective}(\varOptimal)\}$, since the expected-value part has
singleton subdifferential (\cref{proposition:prop_functional_expected_value}),
the variance part is differentiable
(\cref{corollary:cor_functional_variance}) with a gradient that is
continuous by its explicit formula {(convergence in Wasserstein in the argument
and in the value)}, and the sum rule
(\cref{proposition:prop_calculus}) applies.}
Then, by~\cref{theorem:th_first_order_optimality_conditions}, at any optimal $\varOptimal \in \Pp{2}{\reals^d}$ an optimal transport plan $\varPlanGeneric \in \setPlansOptimal{\varOptimal}{\varReference}$, a Lagrange multiplier $\lambda \geq 0$, and a ``sum coupling'' $\varTangentCoupling \in \setPlansCommonMarginal{\gradient{}{\objective}(\varOptimal)}{2\lambda\pushforward{(\projection{1}{}, \projection{1}{} - \projection{2}{})}{\varPlanGeneric}}{1}$
exist so that $\localZero{\varOptimal} \localEquiv{\varOptimal} \gradient{}{\objective}(\varOptimal) \localSum{\varTangentCoupling}2\lambda\pushforward{(\projection{1}{}, \projection{1}{} - \projection{2}{})}{\varPlanGeneric}$. Equivalently,
\begin{align*}
    0 &=
    \min_{\varTangentCoupling \in \setPlansCommonMarginal{\gradient{}{\objective}(\varOptimal)}{2\lambda\pushforward{(\projection{1}{}, \projection{1}{} - \projection{2}{})}{\varPlanGeneric}}{1}}
    \norm{\gradient{}{\objective}(\varOptimal)
    \localSum{\varTangentCoupling}
    2\lambda\pushforward{(\projection{1}{}, \projection{1}{} - \projection{2}{})}{\varPlanGeneric}}^2_{\varOptimal}
    \\
    &=
    \min_{\varTangentCoupling \in \setPlansCommonMarginal{\gradient{}{\objective}(\varOptimal)}{2\lambda\pushforward{(\projection{1}{}, \projection{1}{} - \projection{2}{})}{\varPlanGeneric}}{1}}
    \int_{\reals^d\times\reals^d}\norm{\varTangentReals}^2\d(\gradient{}{\objective}(\varOptimal)
    \localSum{\varTangentCoupling}
    2\lambda\pushforward{(\projection{1}{}, \projection{1}{} - \projection{2}{})}{\varPlanGeneric})(\bar x, \varTangentReals)
    \\
    &=
    \min_{\varTangentCoupling \in \setPlansCommonMarginal{\gradient{}{\objective}(\varOptimal)}{2\lambda\pushforward{(\projection{1}{}, \projection{1}{} - \projection{2}{})}{\varPlanGeneric}}{1}}
    \int_{\reals^d\times\reals^d\times\reals^d}\norm{\varTangentReals_1 + \varTangentReals_2}^2\d\varTangentCoupling(x, \varTangentReals_1, \varTangentReals_2).
\end{align*}
Let now $\varTangentCoupling$ be the minimizer. Then,
\begin{align*}
    0 &=
    \int_{\reals^d\times\reals^d\times\reals^d}\norm{-(1 + \rho(\innerProduct{\theta}{x} - \expectedValue{z \sim \varOptimal}{\innerProduct{\theta}{z}}))\theta + \varTangentReals_2}^2\d\varTangentCoupling(x, \varTangentReals_1, \varTangentReals_2)
    \\
    &=
    \int_{\reals^d\times\reals^d}\norm{-(1 + \rho(\innerProduct{\theta}{x} - \expectedValue{z \sim \varOptimal}{\innerProduct{\theta}{z}}))\theta + \varTangentReals_2}^2\d(\pushforward{\projection{13}{}}{\varTangentCoupling})(x, \varTangentReals_2)
    \\
    &=
    \int_{\reals^d\times\reals^d}\norm{-(1 + \rho(\innerProduct{\theta}{x} - \expectedValue{z \sim \varOptimal}{\innerProduct{\theta}{z}}))\theta + 2\lambda\varTangentReals}^2\d(\pushforward{(\projection{1}{}, \projection{1}{} - \projection{2}{})}{\varPlanGeneric})(x, \varTangentReals)
    \\
    &=
    \int_{\reals^d\times\reals^d}\norm{-(1 + \rho(\innerProduct{\theta}{x} - \expectedValue{z \sim \varOptimal}{\innerProduct{\theta}{z}}))\theta + 2\lambda(x - y)}^2\d\varPlanGeneric(x, y).
\end{align*}
Thus, if $\varOptimal \in \closedWassersteinBall{\varReference}{\varepsilon}{2}$ is optimal, we must have
\begin{equation}
\label{equation:optimality-constrained-simple}
(1 + \rho(\innerProduct{\theta}{x} - \expectedValue{z \sim \varOptimal}{\innerProduct{\theta}{z}}))\theta = 2\lambda(x - y)\quad\almostEverywhere{\varPlanGeneric}
\end{equation}
We now proceed in three steps:

\emph{Step 1.}
We take the expected value with respect to $\varPlanGeneric$ of both sides of \eqref{equation:optimality-constrained-simple} to get
\begin{equation}
\label{eq:theta-eq}
    \theta
    =
    2\lambda(\expectedValue{x \sim \varOptimal}{x} - \expectedValue{y \sim \varReference}{y}).
\end{equation}
If $\lambda = 0$, it must also be $\theta = 0$, so that either there is no solution (if $\theta\neq 0$, consistently with~\cref{example:ex_infimum_not_attained}) or any $\varGeneric \in \Pp{2}{\reals^d}$ is optimal, by direct inspection.
Suppose $\theta\neq0$ and, thus, $\lambda>0$, so that applying the linear map $\innerProduct{\theta}{\cdot}$ to both sides of \eqref{eq:theta-eq} we get:
\begin{equation}
\label{equation:optimality-constrained-simple:expected-value}
    \expectedValue{x \sim \varOptimal}{\innerProduct{\theta}{x}} = \frac{\norm{\theta}^2}{2\lambda} + \expectedValue{y \sim \varReference}{\innerProduct{\theta}{y}}.
\end{equation}

\emph{Step 2.}
We square both sides of \eqref{equation:optimality-constrained-simple} and take the expected value with respect to $\varPlanGeneric$. The right-hand side gives
    $\expectedValue{(x, y)\sim \varPlanGeneric}{\norm{2\lambda(x-y)}^2}
    =4\lambda^2\expectedValue{(x, y)\sim \varPlanGeneric}{\norm{x-y}^2}
    =4\lambda^2\varepsilon^2$,
    where the last equality comes from $\lambda > 0$ and so $\varOptimal \in \boundary{\wassersteinBall{\varReference}{\varepsilon}{2}}{}$.
    The squared norm of the left-hand side, instead, reads in expectation as
        $\expectedValue{(x, y)\sim \varPlanGeneric}{\norm{(1 + \rho(\innerProduct{\theta}{x} - \expectedValue{z \sim \varOptimal}{\innerProduct{\theta}{z}}))\theta}^2}
        =
        \norm{\theta}^2\left(1 + \rho^2\expectedValue{(x, y)\sim \varPlanGeneric}{\left(\innerProduct{\theta}{x} - \expectedValue{z \sim \varOptimal}{\innerProduct{\theta}{z}}\right)^2}\right)
        =
        \norm{\theta}^2\left(1 + \rho^2\variance{x\sim \varOptimal}{\innerProduct{\theta}{x}}\right)$.
Overall,
\begin{equation}
\label{equation:optimality-constrained-simple:variance}
\variance{x\sim \varOptimal}{\innerProduct{\theta}{x}} = \frac{4\lambda^2\varepsilon^2 - \norm{\theta}^2}{\rho^2\norm{\theta}^2}.
\end{equation}

\emph{Step 3.}
We write \eqref{equation:optimality-constrained-simple} as
$
    (\rho\norm{\theta}^2 - 2\lambda)\innerProduct{\theta}{x} + \text{const.} = -2\lambda\innerProduct{\theta}{y}
$
and evaluate the variance with respect to $\varPlanGeneric$ on both sides to get
\begin{equation}
\label{equation:optimality-constrained-simple:variance:2}
    (\rho\norm{\theta}^2 - 2\lambda)^2\variance{x\sim \varOptimal}{\innerProduct{\theta}{x}} = 4\lambda^2 \variance{y \sim \varReference}{\innerProduct{\theta}{y}}.
\end{equation}
Thus, if $\variance{y \sim \varReference}{\innerProduct{\theta}{y}} > 0$, then $\lambda \neq \frac{\rho\norm{\theta}^2}{2}$, and~\eqref{equation:optimality-constrained-simple:variance} and \eqref{equation:optimality-constrained-simple:variance:2} yield $\lambda$ from a polynomial equation (cf. \cite[Proposition 4.7]{Lanzetti2022}).
In this case, the optimal cost is
$\frac{\norm{\theta}^2}{2\lambda} + \expectedValue{y \sim \varReference}{\innerProduct{\theta}{y}} + \frac{2\rho\lambda^2}{(\rho\norm{\theta}^2 - 2\lambda)^2} \variance{y \sim \varReference}{\innerProduct{\theta}{y}}$.
If, instead, $\variance{y \sim \varReference}{\innerProduct{\theta}{y}} = 0$ (and so $\varReference$ is not absolutely continuous), then \eqref{equation:optimality-constrained-simple:variance:2} yields two cases.
If $\lambda = \frac{\rho\norm{\theta}^2}{2}$, then \eqref{equation:optimality-constrained-simple:expected-value} and \eqref{equation:optimality-constrained-simple:variance} yield the optimal cost
    $\expectedValue{y \sim \varReference}{\innerProduct{\theta}{y}} +
    \frac{\rho\norm{\theta}^2\varepsilon^2}{2} + \frac{1}{2\rho}$.
If $\variance{x\sim \varOptimal}{\innerProduct{\theta}{x}} = 0$, then \eqref{equation:optimality-constrained-simple:variance} gives $\lambda = \frac{\norm{\theta}}{2\varepsilon}$ and \eqref{equation:optimality-constrained-simple:expected-value} yields the optimal cost
    $\expectedValue{y \sim \varReference}{\innerProduct{\theta}{y}} + \norm{\theta}\varepsilon$.
With $\lambda$, the optimal solution follows from inverting the affine relation (cf. \eqref{equation:optimality-constrained-simple}):
{\begin{equation*}
    y = T(x) = {\left(I - \frac{\rho\theta\transpose{\theta}}{2\lambda}\right)}x - {\frac{1}{2\lambda}\left({1 - \rho\expectedValue{y \sim \varReference}{\innerProduct{\theta}{y}} - \frac{\rho\norm{\theta}^2}{2\lambda}}\right)\theta}.
\end{equation*}}
Finally, since $\pushforward{(\identity{}, T)}{\varOptimal}\in \setPlansOptimal{\varOptimal}{\varReference}$, $T$ is an optimal transport map. By~\cite[Theorem 5.10]{villani2009optimal}, it is also a gradient of a convex function, which implies $\lambda \geq \frac{\rho\norm{\theta}^2}{2}$. In particular, our results generalize \cite{Lanzetti2022}.
\end{example}

\section{Applications}\label{sec:applications}
We demonstrate the practical appeal of our first-order optimality conditions by discussing several applications in machine learning and nonlinear \gls*{acr:dro}.
Importantly, our optimality conditions have already contributed to the development of several state-of-the-art algorithms, including learning diffusion processes from population data \cite{terpin2024learning} and energy-based generative modeling \cite{balcerak2025energy}.

\subsection{Machine learning}
In \cite{terpin2024learning}, we showcase how our first-order optimality conditions improve the computational efficiency of learning the dynamics of particles undergoing diffusion, bypassing cumbersome bi-level optimization procedures \cite{bunne2022proximal,terpin2024learning}. In this section, we study applications to drug discovery and distribution fitting.

\subsubsection{Drug discovery}
In Euclidean settings, the proximal operator appears as an often closed-form subroutine of various first-order optimization algorithms \cite{bauschke2011convex}. Recently, it has found applications in learning the dynamics of particles undergoing diffusion \cite{bunne2022proximal,terpin2024learning},
single-cell perturbation analysis \cite{bunne2023single-cell-perturbation}, and
drug discovery~\cite[Section 6]{alvarez2021optimizing}. We use the first-order optimality conditions provided in this work to study the proximal operator, and we showcase how these results can be applied to molecular discovery. More precisely, following~\cite[Section 6]{alvarez2021optimizing}, the task of drug discovery and drug repurposing can be cast as an iterative process that aims to increase the \emph{drug-likeness} $\expectedValue{x \sim \varGeneric}{V(x)}$ of a distribution of molecules $\varGeneric$ while staying close to the original distribution $\varGeneric_t$, which (up to replacing $V$ with $-V$) gives the optimization problem
\[
    \varGeneric_{t+1}
    =
    \argmin_{\varGeneric\in\Pp{2}{\reals^d}}
    \expectedValue{x \sim \varGeneric}{V(x)} + \frac{1}{2}
    \wassersteinDistance{\varGeneric}{\varGeneric_t}{2}^2.
\]
The iterative approach is preferred in practice over a one-shot solution due to the nonconvexity of the drug-likeness potential $V$.
We now study this optimization problem.
Unlike \cite[Section 6]{alvarez2021optimizing}, we do not require any convex proxy for it:
\begin{example}[Proximal operator in the Wasserstein space]\label{ex:proximal_operator}
Let $\varFixed\in\Pp{2}{\reals^d}$ and $V:\reals^d\to\reals$ satisfy the assumptions of \cref{proposition:prop_functional_expected_value}\ref{prop_functional_expected_value_differentiable}.
Consider the proximal operator in the Wasserstein space, defined by
\begin{equation}\label{eq:proximal_operator:definition}
    {\argmin_{\varGeneric\in\Pp{2}{\reals^d}}
    \objective(\varGeneric),
    \quad\text{where}\quad\objective(\varGeneric)
    \coloneqq
    \expectedValue{x\sim\varGeneric}{V(x)}
    +\frac{1}{2}\wassersteinDistance{\varGeneric}{\varFixed}{2}^2.}
\end{equation}
We aim to characterize optimal solutions of~\eqref{eq:proximal_operator:definition}.
Since $V$ is differentiable, the sum rule (\cref{proposition:prop_calculus}\ref{prop_calculus_sum_rule} with $\objective_1(\varGeneric)=\frac{1}{2}\wassersteinDistance{\varGeneric}{\varFixed}{2}^2$ and $\objective_2(\varGeneric)=\expectedValue{x\sim\varGeneric}{V(x)}$) gives
    $\subgradient{\objective}(\varGeneric)
    \subseteq
    \left\{\pushforward{(\identity{},\gradient{}V)}{\varGeneric} \localSum{\varTangentCoupling} \varTangent \st \varTangent \in \pushforward{(\projection{1}{},\projection{1}{}-\projection{2}{})}{\setPlansOptimal{\varGeneric}{\varFixed}}, \varTangentCoupling \in \setPlansCommonMarginal{\pushforward{(\identity{},\gradient{}V)}{\varGeneric}}{\varTangent}{1}\right\}$.
By~\cref{theorem:th_fermat_rule}, the subgradient of $\objective{}$ vanishes at an optimal solution $\varOptimal$. Thus, at optimality there exist $\varPlanGeneric\in\setPlansOptimal{\varOptimal}{\varFixed}$ and $\varTangentCoupling \in \setPlansCommonMarginal{\pushforward{(\identity{},\gradient{}V)}{\varOptimal}}{\pushforward{(\projection{1}{}, \projection{1}{} - \projection{2}{})}{\varPlanGeneric}}{1}$
so that $\localZero{\varOptimal} \localEquiv{\varOptimal}\pushforward{(\identity{},\gradient{}V)}{\varOptimal}
\localSum{\varTangentCoupling}
\pushforward{(\projection{1}{},\projection{1}{}-\projection{2}{})}{\varPlanGeneric}$.
After the (by now) usual steps, we conclude that this equivalently means that
for almost all $(x,y)\in\supp(\varPlanGeneric)$, it necessarily holds
    $0=\gradient{}V(x)+(x-y)=\gradient{}{}\left(V(x)+\frac{\norm{x-y}^2}{2} \right)$,
which recovers the standard optimality condition for the proximal operator in Euclidean spaces.
By direct inspection, we conclude that the cost-minimizing choice is $x\in\argmin_{z\in\reals^d} V(z)+\frac{\norm{z-y}^2}{2}$ for $\almostEverywhere{\varPlanGeneric}$ $y$,
    $\objective(\varOptimal)
    =
    \expectedValue{y\sim\varFixed}{\min_{x\in\reals^d}V(x)+\frac{1}{2}\norm{x-y}^2}$.
In particular, the proximal operator~\eqref{eq:proximal_operator:definition} in $\Pp{2}{\reals^d}$ can be expressed in terms of the proximal operator on $\reals^d$. In the special case where the function $x\mapsto V(x)+\frac{\norm{x}^2}{2}$ is strictly convex, we conclude $\varOptimal=\pushforward{(\identity{}+\gradient{}V)^{-1}}{\varFixed}$.
\end{example}

Importantly, \cref{ex:proximal_operator} suggests a closed-form solution that can be incorporated into, e.g., \cite{alvarez2021optimizing} to improve the efficiency of the numerical procedure. The numerical benefits of using reformulations inspired by our first-order optimality conditions have been explored in detail in, e.g., \cite{terpin2024learning}.

\subsubsection{Learning Gaussians}
\label{example:learning-gaussian-mm}
\begin{figure}
    \centering
    \begin{minipage}{.32\textwidth}
    \centering
    \includegraphics[width=\linewidth]{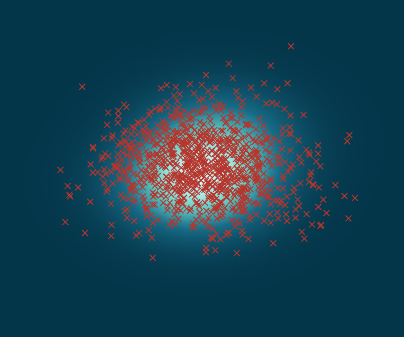}
    \end{minipage}
    \hfill
    \begin{minipage}{.32\textwidth}
    \centering
    \includegraphics[width=\linewidth]{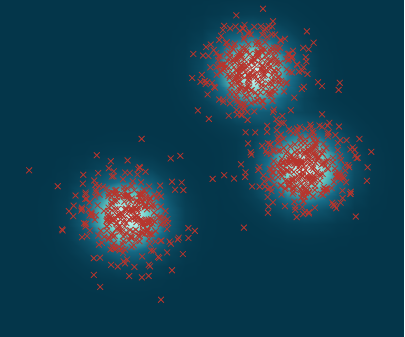}
    \end{minipage}
    \hfill
    \begin{minipage}{.32\textwidth}
    \centering
    \includegraphics[width=\linewidth]{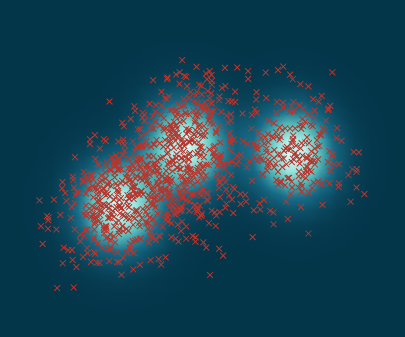}
    \end{minipage}
    \caption{Data (red crosses) generated from mixtures with one, three, and five Gaussian components (left to right) are fitted using the method outlined in \Cref{example:learning-gaussian-mm} with three components.}
    \label{fig:mixture-gaussian}
\end{figure}

Gaussian mixture models are among the most widespread techniques for distribution fitting \cite{reynolds2009gaussian,terpin2024learning,yan2023learning}. It is well known (e.g., see~\cite{yan2023learning} and references therein) that this problem can be formulated as an optimization problem in the space of probability measures:
Since a Gaussian mixture model can be written as $\varGeneric\ast\phi$ with $\phi(x)$ being the Gaussian kernel\footnote{Here, we consider the mixture of Gaussians with unit variance, but the discussion in this section easily extends to general mixtures.}
$ \frac{1}{(2\pi)^{d/2}}\exp(-\frac{\norm{x}^2}{2})$, the learning problem can be cast as the optimization problem
	$\min_{\varGeneric\in\Pp{2}{\reals^d}}
	-\sum_{i=1}^n\log\left(\int_{\reals^d}\phi(x-X_i)\d\varGeneric(x)\right)$,
where $\{X_1,\ldots,X_n\}$ is a given dataset.
Since the objective is differentiable, \cref{theorem:th_fermat_rule} yields the necessary optimality condition
\begin{equation}
	0
	=
    \int_{\reals^d}
    \norm{
	-\sum_{i=1}^n\frac{\gradient{}\phi(x-X_i)}{\int_{\reals^d}\phi(y-X_i)\d\varOptimal(y)}
 }^2
    \d\varOptimal(x)
	=
 \int_{\reals^d}
 \norm{
	\sum_{i=1}^n\frac{(x-X_i)\exp(-\frac{\norm{x-X_i}^2}{2})}{\int_{\reals^d}\exp(-\frac{\norm{y-X_i}^2}{2})\d\varOptimal(y)}
}^2
 \d\varOptimal(x),\label{eq:opt-condition}
\end{equation}
where we used the sum and chain rules (\cref{proposition:prop_calculus}) and differentiability of $\varGeneric \mapsto \int_{\reals^d}\phi(x-X_i)\d\varGeneric(x)$ (\cref{proposition:prop_functional_expected_value}\ref{prop_functional_expected_value_differentiable}) to evaluate the subgradient.
At this point, one can parametrize $\varGeneric$ with finitely many samples $x_1,\ldots,x_m$, with weights $\rho_1,\ldots,\rho_m$, and search for the $x_1,\ldots,x_m$ and $\rho_1,\ldots,\rho_m$ which best ``fit'' the first-order optimality condition. We showcase the results of the method in \cref{fig:mixture-gaussian}. Most importantly, our method also provides a simple and computationally efficient way to numerically assess the output of other distribution fitting methods: For approximately optimal solutions, evaluating \eqref{eq:opt-condition} must result in negligible values.

\subsection{Nonlinear mean-variance optimization}

In this section, we extend the existing duality result for linear mean-variance distributionally robust optimization \cite{blanchet2019quantifying,blanchet2022optimal,kuhn2019wasserstein,yue2022linear,gao2023distributionally} to the nonlinear case.
Specifically, for $\varReference \in \Pp{2}{\reals^d}, \rho>0, \varepsilon > 0$, and a differentiable function $V: \reals^d \to \reals$ satisfying the assumptions of \cref{corollary:cor_functional_variance}, consider the problem of evaluating the \emph{worst-case} risk over a Wasserstein ball:
\begin{equation}
\label{equation:dro-meanvar:primal}
{%
    \objective(\varOptimal)
    =
    \inf_{\varGeneric \in \closedWassersteinBall{\varReference}{\varepsilon}{2}}
    \objective(\varGeneric)
    \quad\text{where}\quad\objective(\varGeneric)
    \coloneqq
    -\left(\expectedValue{x \sim \varGeneric}{V(x)} + \rho\variance{x \sim \varGeneric}{V(x)}\right).}
\end{equation}
We use our results to show that any candidate optimal solution $\varOptimal \in \Pp{2}{\reals^d}$ to \eqref{equation:dro-meanvar:primal} yields the worst-case cost
\begin{equation}
\label{equation:dro-meanvar:dual}
\objective(\varOptimal) = -\inf_{\substack{\lambda \geq 0\\\beta \in \reals}}\; \lambda \varepsilon^2 + \int_{\reals^d} \max_{y \in \reals^d} \left\{ V(y) + \rho(V(y) - \beta)^2 - \lambda \norm{y - \hat{x}}^2\right\}\d\varReference(\hat{x}).
\end{equation}
Importantly, \eqref{equation:dro-meanvar:dual} reduces the infinite-dimensional problem in the space of probability measures \eqref{equation:dro-meanvar:primal} to a two-dimensional problem over the reals.
Moreover, with $\lambda^\ast, \beta^\ast$ being the minimizers above, there exists $\varPlanGeneric \in \setPlansOptimal{\varOptimal}{\varReference}$ such that
\begin{equation}
\label{equation:dro-meanvar:selection}
x^\ast \in \argmax_{y \in \reals^d} \left\{ V(y) + \rho(V(y) - \beta^\ast)^2 - \lambda^\ast \norm{y - \hat{x}}^2\right\}\quad \forall (x^\ast, \hat{x}) \in \supp(\varPlanGeneric).
\end{equation}
Thus, the optimizers $\lambda^\ast,\beta^\ast$ can be used to construct an optimal solution of \eqref{equation:dro-meanvar:primal} by allocating the mass at $\hat x$ $\almostEverywhere{\varReference}$ to the corresponding set of maximizers in \eqref{equation:dro-meanvar:selection}.

To prove~\eqref{equation:dro-meanvar:dual}, we establish the two inequalities separately.
For ``$\leq$'', \Cref{theorem:th_first_order_optimality_conditions} combined with \cref{proposition:prop_functional_expected_value,proposition:prop_calculus,corollary:cor_functional_variance,proposition:cor_subgradients_horizon,proposition:prop_normal_cone_wasserstein_ball} (in particular, the constraint qualification holds since each term of $\objective$ is Lipschitz on bounded sets and the horizon subgradient is thus trivial, cf. \cref{proposition:cor_subgradients_horizon,proposition:prop_subgradients_horizon_trivial}, and the closed Wasserstein ball is \gls*{acr:snc} by \cref{proposition:prop_normal_cone_wasserstein_ball}) requires that, for any optimal solution $\varOptimal\in \Pp{2}{\reals^d}$, there are $\varPlanGeneric \in {\setPlansOptimal{{\varOptimal}}{\varReference}}$ and $\varTangent = 2\lambda^\ast\pushforward{(\projection{1}{}, \projection{1}{} - \projection{2}{})}{\varPlanGeneric}$ with $\lambda^\ast > 0$ (in particular, $\wassersteinDistance{\varOptimal}{\varReference}{2} = \varepsilon$ and $\varTangent\in\normalConeStrong{\closedWassersteinBall{\varReference}{\varepsilon}{2}}{\varOptimal}$; for $\lambda^\ast = 0$, see \cref{example:ex_infimum_not_attained}) such that
{
    $\localZero{{\varOptimal}}
    \in
    \pushforward{(\projection{1}{}, -\comp{\projection{1}{}}{\gradient{}V}(1 + 2\rho(\comp{\projection{1}{}}{V} - \expectedValue{x\sim\varOptimal}{V(x)})) + \projection{2}{})}{\varTangent}
    $
    or, equivalently,
    $\localZero{{\varOptimal}}
    \in
    \pushforward{(\projection{1}{}, -\comp{\projection{1}{}}{\gradient{}V}(1 + 2\rho(\comp{\projection{1}{}}{V} - \expectedValue{x\sim\varOptimal}{V(x)})) + 2\lambda^\ast(\projection{1}{}-\projection{2}{}))}{\varPlanGeneric}
    $}.
This implies that for all $(x^\ast, \hat{x}) \in \supp(\varPlanGeneric)$ we have
$
0=- \gradient{}{V}(x^\ast)(1 + 2\rho(V(x^\ast) - \expectedValue{x\sim\varOptimal}{V(x)}))+2\lambda^\ast(x^\ast - \hat{x})$
{, which, following~\cref{ex:proximal_operator}, allows us to conclude
$x^\ast \in \argmin_{y \in \reals^d} \left\{
-V(y) - \rho(V(y) - \expectedValue{x\sim\varOptimal}{V(x)})^2
+ \lambda^\ast\norm{y - \hat{x}}^2
\right\}$.
Overall, we have
\begin{equation}\label{eq:mean-variance:local_optimality_condition}
	x^\ast \in \argmax_{y \in \reals^d} \left\{
	V(y) + \rho(V(y) - \expectedValue{x\sim\varOptimal}{V(x)})^2
	- \lambda^\ast\norm{y - \hat{x}}^2
	\right\}
	\quad \forall (x^\ast, \hat{x}) \in \supp(\varPlanGeneric).
\end{equation}}
Thus, with $\variance{x\sim\varGeneric}{V(x)}=\expectedValue{x\sim\varGeneric}{(V(x)-\expectedValue{y\sim\varGeneric}{V(y)})^2}$, we have
\begin{align*}
    \objective(\varOptimal)
    &=
    -\expectedValue{x^\ast \sim \varOptimal}{V(x^\ast)+\rho\left(V(x^\ast)-\expectedValue{x\sim\varOptimal}{V(x)}\right)^2}
    \\
    \overset{\eqref{eq:mean-variance:local_optimality_condition}}&{=}
    \begin{aligned}
    -\int_{\reals^d\times\reals^d}
    V(x^\ast)+&\rho\left(V(x^\ast)-\expectedValue{x\sim\varOptimal}{V(x)}\right)^2
    - \lambda^\ast\norm{x^\ast - \hat{x}}^2 + \lambda^\ast\norm{x^\ast - \hat{x}}^2
    \d\varPlanGeneric(x^\ast,\hat x)
    \end{aligned}
    \\
    &=
    -\lambda^\ast\varepsilon^2
    -\int_{\reals^d}\max_{y \in \reals^d}
    \left\{
    V(y) + \rho\left(V(y) - \expectedValue{x\sim\varOptimal}{V(x)}\right)^2
    - \lambda^\ast\norm{y - \hat{x}}^2
    \right\}
    \d\varReference(\hat{x})
    \\
    &\leq
    \sup_{\substack{\lambda\geq 0 \\ \beta\in\reals}}
    -\lambda\varepsilon^2
    -\int_{\reals^d}\max_{y \in \reals^d}
    \left\{
    V(y) + \rho\left(V(y) - \beta\right)^2
    - \lambda\norm{y - \hat{x}}^2
    \right\}
    \d\varReference(\hat{x}),
\end{align*}
which establishes the first inequality.
For ``$\geq$'', we rewrite $\variance{x\sim\varGeneric}{V(x)}$ as the minimization $\inf_{\beta\in\reals}\expectedValue{x\sim\varGeneric}{(V(x)-\beta)^2}$ and use a standard minimax argument:
\begin{align*}
	\objective(\varOptimal)
	&=
	\inf_{\varGeneric\in\Pp{2}{\reals^d}}\sup_{\lambda\geq 0} -\expectedValue{x\sim\varGeneric}{V(x)}-\rho\inf_{\beta\in\reals}\expectedValue{x\sim\varGeneric}{(V(x)-\beta)^2}
		+\lambda(\wassersteinDistance{\varGeneric}{\varReference}{2}^2-\varepsilon^2)
	\\
	&=
	\inf_{\varGeneric\in\Pp{2}{\reals^d}}\sup_{\substack{\lambda\geq 0 \\ \beta\in\reals}} -\expectedValue{x\sim\varGeneric}{V(x)+\rho(V(x)-\beta)^2}+\lambda(\wassersteinDistance{\varGeneric}{\varReference}{2}^2-\varepsilon^2)
	\\
	&\geq
	\sup_{\substack{\lambda\geq 0 \\ \beta\in\reals}}
	\inf_{\varGeneric\in\Pp{2}{\reals^d}} -\expectedValue{x\sim\varGeneric}{V(x)+\rho(V(x)-\beta)^2}+\lambda(\wassersteinDistance{\varGeneric}{\varReference}{2}^2-\varepsilon^2)
	\\
	&=
	\sup_{\substack{\lambda\geq 0 \\ \beta\in\reals}}
	-\lambda\varepsilon^2
    -\int_{\reals^d}\max_{y \in \reals^d}
    \left\{
    V(y) + \rho\left(V(y) - \beta\right)^2
    - \lambda\norm{y - \hat{x}}^2
    \right\}
    \d\varReference(\hat{x}),
\end{align*}
where the last equality follows from the proximal operator in~\cref{ex:proximal_operator}. This establishes the second inequality and concludes the proof of~\eqref{equation:dro-meanvar:dual}.

\bibliography{references}

\newpage
\appendix
\crefname{appendix}{Appendix}{Appendices}
\Crefname{appendix}{Appendix}{Appendices}
\crefalias{section}{appendix}
\crefalias{subsection}{appendix}
\crefalias{subsubsection}{appendix}
\section{Auxiliary technical results}\label{app:technical-results}
\subsection{Basics of calculus and optimal transport}
This section gathers auxiliary technical results of limited novelty used throughout the paper.
We start with some properties of functionals on $\reals^d$:
\includeElement{proposition}{prop_subgradient_convex}{true}{true}

The second result extends \cite[Lemma 6.1]{terpin2023dynamic} to the case where one marginal is common. In particular, the proof is analogous, and the same proof allows us to state the result ``pointwise'', rather than just for the infimum. The extension to multiple marginals, or to several common marginals, is straightforward, so we omit the proof.

\includeElement{lemma}{lemma_pushforward_and_optimal_transport}{true}{false}

{
\begin{remark}[Extension to countably many spaces]
\label{remark:rem_lemma_pushforward_and_optimal_transport}
    The proof of \cite[Lemma 6.1]{terpin2023dynamic} extends to countably many spaces, provided that the limit distribution arising from the iterated application of \cite[Gluing Lemma]{villani2009optimal} is well defined. For this, it suffices to invoke the Kolmogorov extension theorem \cite[Section 2.4]{tao2011introduction}.
\end{remark}
}

{The next result concerns the stability of the set of
transport plans:}

\begin{lemma}\label{lemma:set_plans}
  Fix $m\in\naturals$, $m\geq2$, and $\varOther_1,\dots,\varOther_{m-1}\in\Pp{2}{\reals^d}$ and {consider the
  set-valued map}
  $\varGeneric\mapsto G(\varGeneric)=\setPlans{\varGeneric}{\varOther_1,\dots,\varOther_{m-1}}\subseteq\Pp{2}{(\reals^d)^m}$,
  defined on $\Pp{2}{\reals^d}$.
    Then, $G$ is lower hemicontinuous in Wasserstein.
\end{lemma}

\begin{proof}{Proof.}
    We prove the case with two marginals ($m=2$), the general case being analogous.
    For lower hemicontinuity, let $\varGeneric_n\convergenceWasserstein\varGeneric$ and $\varPlanGeneric\in G(\varGeneric)=\setPlans{\varGeneric}{\varOther_1}$. We need to show that there is a sequence $\varPlanGeneric_{n}\in G(\varGeneric_{n})=\setPlans{\varGeneric_n}{\varOther_1}$ so that $\varPlanGeneric_{n}\convergenceWasserstein\varPlanGeneric$.
    To start, since $\varGeneric_n\convergenceWasserstein\varGeneric$, there is a sequence $\bar\varPlanGeneric_n\in\setPlans{\varGeneric}{\varGeneric_n}$ so that $\int_{\reals^d\times\reals^d}\norm{x-x_n}^2\d\bar\varPlanGeneric_n(x,x_n) \to 0$.
    Via \cite[Gluing Lemma]{villani2009optimal}, define $\varTangentCoupling_n\in\setPlansCommonMarginal{\bar\varPlanGeneric_n}{\varPlanGeneric}{1}$ by ``gluing'' $\bar\varPlanGeneric_n$ and $\varPlanGeneric$; in particular, $\pushforward{(\projection{1}{},\projection{2}{})}{\varTangentCoupling_n}=\bar\varPlanGeneric_n$ and $\pushforward{(\projection{1}{},\projection{3}{})}{\varTangentCoupling_n}=\varPlanGeneric$.
    Define $\varPlanGeneric_n=\pushforward{(\projection{2}{},\projection{3}{})}{\varTangentCoupling_n}$. By construction, $\varPlanGeneric_n\in G(\varGeneric_n)=\setPlans{\varGeneric_n}{\varOther_1}$. We now claim $\varPlanGeneric_n\convergenceWasserstein\varPlanGeneric$. With the candidate transport plan $\beta_n=\pushforward{(\projection{2}{},\projection{3}{},\projection{1}{},\projection{3}{})}{\varTangentCoupling_n}\in\setPlans{\varPlanGeneric_n}{\varPlanGeneric}$, we have
    \begin{align*}
        \wassersteinDistance{\varPlanGeneric_n}{\varPlanGeneric}{2}^2
        &\leq
        \int_{\reals^d\times\reals^d\times\reals^d\times\reals^d}
        \norm{x_n-x}^2+\norm{y_n-y}^2\d\beta_n(x_n,y_n,x,y)
        \\
        &=
        {\int_{\reals^d\times\reals^d\times\reals^d}
        \norm{x_n-x}^2+\norm{y-y}^2\d\varTangentCoupling_n(x,x_n,y)}
        \\
        &=
        {\int_{\reals^d\times\reals^d\times\reals^d}
        \norm{x_n-x}^2\d\varTangentCoupling_n(x,x_n,y)}
        \\
        &=
        \int_{\reals^d\times\reals^d}
        \norm{x_n-x}^2\d\bar\varPlanGeneric_n(x,x_n)
        \to 0.
    \end{align*}
    Thus, $\varPlanGeneric_n\convergenceWasserstein\varPlanGeneric$, which concludes the proof.
\end{proof}

\subsection{Geometry of the Wasserstein space}
We recall from \cite{gigli2008geometry} several properties of the operations in \cref{definition:def_almost_linear_structure_wass}:
\begin{proposition}[Sum and scale]
\label{proposition:prop_sum_scale_well_defined}
{Let $\varTangent_1, \varTangent_2 \in \tangentCone{\Pp{2}{\reals^d}}{\varFixed}$. Then,
\begin{enumerate}
    \item For all $\tau \in \reals$, $\tau\varTangent_1\in\tangentCone{\Pp{2}{\reals^d}}{\varFixed}$.
    \item For $\varTangentCoupling\in\setPlansCommonMarginal{\varTangent_1}{\varTangent_2}{1}$,
$
\varTangent_1
\localSum{\varTangentCoupling}
\varTangent_2
\coloneqq
\pushforward{(\projection{1}{}, \projection{2}{} + \projection{3}{})}{\varTangentCoupling}
\in\tangentCone{\Pp{2}{\reals^d}}{\varFixed}$.
\end{enumerate}}
\end{proposition}

\begin{proof}{Proof of~\cref{proposition:prop_sum_scale_well_defined}.}
  See \cite[Proposition 4.25]{gigli2008geometry} for (i) and \cite[Proposition 4.29]{gigli2008geometry} and the preceding discussion for (ii).
\end{proof}

Next, we study the geometry of the Wasserstein space:

\includeElement{proposition}{prop_wasserstein_geom}{true}{true}

{Finally, we collect how a \emph{deterministic coordinate}---a
Dirac marginal in the last coordinate of the base measure, the situation
arising systematically for epigraphs (cf.
\cref{section:variational-geometry:epigraphs})---interacts with the
Wasserstein geometry:}

\includeElement{lemma}{lemma_deterministic_coordinate}{true}{true}

{In words, a deterministic coordinate is ``transparent'' to the
Wasserstein geometry: it contributes an additive, transport-plan-independent term to
transport costs, it neither creates nor destroys tangency, and it enters
norms and inner products separately from the other coordinates, through the
deterministic velocity $c$ and the expectation of the last velocity component.}

\section{Proofs for \texorpdfstring{\cref{section:wasserstein_space}}{Section 2}}\label{appendix:sec2:proofs}

\includeElement{proposition}{prop_differentiability_properties}{duplicateStatementInAppendix}{proofInAppendix}
\includeElement{proposition}{prop_wassersteindistance_regular_subgradient_must_be_plan}{duplicateStatementInAppendix}{proofInAppendix}
\includeElement{proposition}{prop_wassersteindistance_general_subgradient}{duplicateStatementInAppendix}{proofInAppendix}
\includeElement{proposition}{prop_functional_optimal_transport}{duplicateStatementInAppendix}{proofInAppendix}
\includeElement{proposition}{prop_functional_expected_value}{duplicateStatementInAppendix}{proofInAppendix}
\includeElement{proposition}{prop_functional_interaction}{duplicateStatementInAppendix}{proofInAppendix}
\includeElement{proposition}{prop_calculus}{duplicateStatementInAppendix}{proofInAppendix}
\includeElement{corollary}{cor_functional_variance}{duplicateStatementInAppendix}{proofInAppendix}
\includeElement{proposition}{prop_normal_tangent_cones}{duplicateStatementInAppendix}{proofInAppendix}
\includeElement{proposition}{prop_normal_cones_trivial}{duplicateStatementInAppendix}{proofInAppendix}
\includeElement{lemma}{strong_vs_regular_normal_cone}{duplicateStatementInAppendix}{proofInAppendix}
\includeElement{proposition}{prop_epigraph_normal_tangent}{duplicateStatementInAppendix}{proofInAppendix}
\includeElement{proposition}{prop_subgradients_epigraphical_normals}{duplicateStatementInAppendix}{proofInAppendix}
\includeElement{proposition}{prop_subgradients_horizon_trivial}{duplicateStatementInAppendix}{proofInAppendix}
\includeElement{proposition}{cor_subgradients_horizon}{duplicateStatementInAppendix}{proofInAppendix}
\includeElement{proposition}{th_subgradients_epigraphical_normals}{duplicateStatementInAppendix}{proofInAppendix}
\includeElement{proposition}{prop_normal_cones_constraints_simple}{duplicateStatementInAppendix}{proofInAppendix}
\includeElement{proposition}{prop_normal_cone_wasserstein_ball}{duplicateStatementInAppendix}{proofInAppendix}
\section{Proofs for \texorpdfstring{\cref{section:optimization}}{Section 3}}

Here, we follow the strategy outlined in \cref{section:variational-geometry:epigraphs}, so we start by proving the first-order optimality conditions for a differentiable functional:
\includeElement{proposition}{cor_first_order_optimality_conditions_diff}{true}{proofInAppendix}
{%
The proof of our main result lifts the problem to the epigraph, where the objective becomes a linear functional; the next lemma shows that such functionals are strongly differentiable, so that \cref{proposition:cor_first_order_optimality_conditions_diff} applies.
\includeElement{lemma}{linear_strongly_differentiable}{true}{proofInAppendix}
}

{%
We record a companion to \cref{lemma:linear_strongly_differentiable}: the normal
cone and \gls*{acr:snc} property of a constraint set defined as the level set of an expected value of a linear function. We will prove this result from the definition and use it to prove~\cref{proposition:prop_normal_cones_constraints_simple}.}
\includeElement{lemma}{normal_cone_linear_constraint}{true}{proofInAppendix}

Armed with \cref{proposition:cor_first_order_optimality_conditions_diff} and \cref{lemma:linear_strongly_differentiable}, we are ready to prove our main result and its differentiable case:
\includeElement{theorem}{th_first_order_optimality_conditions}{duplicateStatementInAppendix}{proofInAppendix}
\includeElement{theorem}{th_fermat_rule}{duplicateStatementInAppendix}{proofInAppendix}

\section{Extremal principle and intersection rule}
\label{app:intersection}

In this section, we establish the intersection rule for the strong normal cone: Under an appropriate constraint qualification, we have
$
\normalConeStrong{\Omega_1\cap\Omega_2}{\varFixed} \subseteq \normalConeStrong{\Omega_1}{\varFixed}
+
\normalConeStrong{\Omega_2}{\varFixed}
$.
This intersection rule plays a fundamental role in variational analysis and, among others, in the proof of our first-order optimality conditions in~\cref{theorem:th_first_order_optimality_conditions}.
In the Euclidean case and when all sets are convex, the intersection rule follows from a convex separation principle, which however fails to apply to the infinite-dimensional and nonconvex case presented here.
We therefore follow the strategy of variational analysis in linear spaces (in particular, \WikipediaLink{https://en.wikipedia.org/wiki/Banach_space}{Banach spaces} and \WikipediaLink{https://en.wikipedia.org/wiki/Asplund_space}{Asplund spaces})~\cite{mordukhovich2006variationalI,mordukhovich2006variationalII} and prove the intersection rule via an extremal principle.
An extremal principle provides necessary conditions for local extremal points of general (nonconvex) sets in terms of the normal cone at those points. It can be interpreted as the ``variational version'' of the convex separation principle; see the discussion in~\cite[Section 2]{mordukhovich2006variationalI}.
Importantly, the Wasserstein space is not a linear space (as opposed to Banach and Asplund spaces), but is a metric space, so the results~\cite{mordukhovich2006variationalI,mordukhovich2006variationalII} do not directly apply.
Accordingly, in~\cref{subsec:extremal_principle}, we establish an extremal principle in the Wasserstein space, which we believe to be of independent interest.
Then, we establish the product rule in~\cref{subsec:product_rule} and the intersection rule in~\cref{subsec:intersection_rule}.

\subsection{Extremal principle}\label{subsec:extremal_principle}

We start with a definition of the extremal principle in the Wasserstein space:

\includeElement{definition}{def_extremal}{true}{false}

Informally, local extremality means that two sets of probability measures can be separated by ``moving'', via a pushforward operation, one of them.
In the case of delta distributions on $\reals^d$, we recover the definition of extremal system in $\reals^d$; e.g., see~\cite{mordukhovich2018variational} or~\cite[Definition 2.1]{mordukhovich2006variationalI} in the finite-dimensional case.
We next turn our focus to necessary conditions for $(\Omega_1,\Omega_2,\varFixed)$ to be extremal.

\includeElement{theorem}{th_extremal_principle}{proofInAppendix}{proofInMain}

Informally, \cref{theorem:th_extremal_principle} shows that if a probability measure $\varFixed$ is a locally extremal point of $\Omega_1$ and $\Omega_2$ then at probability measures arbitrarily close to $\varFixed$ there are unit tangent vectors, arbitrarily close to the strong regular normal cones of $\Omega_1$ and of $\Omega_2$, which oppose each other.
The proof of~\cref{theorem:th_extremal_principle} relies on the following technical lemma:

\includeElement{lemma}{lemma_fuzzy_optimality_condition}{proofInAppendix}{proofInAppendix}

We are now ready to prove~\cref{theorem:th_extremal_principle}.
The proof follows the high-level strategy in~\cite[Section 2.5.3]{mordukhovich2006variationalI}, with several important modifications. For instance,
\begin{enumerate}[label=\text{\faBolt}~\arabic*.]
    \item The properties of the Wasserstein distance force us to use the Borwein--Preiss variational principle instead of the Ekeland variational principle (see~\cite[Theorem 2.26]{mordukhovich2006variationalI}).
    \item The tangent vectors are localized, which renders the ones at different points incomparable.
    \item We cannot resort to the linear structure of the space and the properties of norms.
\end{enumerate}

\includeElement{theorem}{th_extremal_principle}{duplicateStatementInAppendix}{proofInAppendix}

{%
\subsection{Product rule}\label{subsec:product_rule}

\includeElement{proposition}{prop_normal_cone_product}{proofInAppendix}{proofInAppendix}

{%
\begin{remark}[Regular vs.\ strong normal cone]
In \cref{proposition:prop_normal_cone_product}\ref{item:product:cylinder} equality holds only for the strong normal cones; for the (regular) normal cone only the inclusion holds.
Indeed, the reverse inclusion tests against variations from \emph{optimal} transport plans, whose $\reals^d$-marginal need not be an optimal transport plan of $\varFixed_2\coloneqq\pushforward{\projection{2}{}}{\varFixed}$, so $\pushforward{\projection{24}{}}{\bar\varTangent}\in\normalConeRegular{\Omega_2}{\varFixed_2}$ does not suffice.
The strong normal cones test against all transport plans, removing the gap.
\end{remark}%
}

\subsection{Intersection rule}\label{subsec:intersection_rule}

We now turn our attention to the intersection rule.

Our first lemma converts an \emph{approximate} strong regular normal, a variation satisfying the strong regular normal cone inequality up to an additive~$\varepsilon$, into an \emph{exact} strong regular normal at a nearby base point.

\newcommand{\epigraphSpace}{\mathbb{X}}
\includeElement{lemma}{eps_to_exact_normal}{true}{true}

Our second lemma provides us with a fuzzy intersection rule, which will be key in establishing the general intersection rule below.

\includeElement{lemma}{fuzzy-intersection}{true}{true}

We are now ready to state and prove the intersection rule for the strong normal cone.

\includeElement{theorem}{intersection_rule_strong}{true}{true}
}

\end{document}